\documentclass[12pt,a4paper,nosumlimits]{amsart}
\usepackage{latexsym}

\calclayout

\global\let\le\undefined
\DeclareMathSymbol{\le}{\mathrel}{AMSa}{"36}     
\global\let\ge\undefined
\DeclareMathSymbol{\ge}{\mathrel}{AMSa}{"3E}     
\global\let\k\undefined
\DeclareMathSymbol{\k}{\mathord}{AMSb}{"7C}      

\newcommand{\ds}{\dots}
\newcommand{\dsm}{\dotsm}
\newcommand{\dsb}{\dotsb}
\newcommand{\cd}{\cdot}
\def\ot{\DOTSB\otimes}
\newcommand{\op}{\oplus}
\def\rarrow{\DOTSB\longrightarrow}
\def\larrow{\DOTSB\longleftarrow}

\def\lrarrow{\DOTSB\,\relbar\joinrel\relbar\joinrel\rightarrow\,}
\def\llarrow{\DOTSB\,\leftarrow\joinrel\relbar\joinrel\relbar\,}

\renewcommand{\:}{\colon}
\newcommand{\oa}{\odot}
\newcommand{\+}{\nobreakdash-}
\newcommand{\1}{\prime}
\renewcommand{\.}{\mskip .5\thinmuskip}

\newcommand{\nbk}{\nobreak}

\newcommand{\ind}{}
\newenvironment{thm}[1]{\medskip\ind{\sf #1.}\sl}{\medskip}

\newenvironment{rem}[1]{\medskip\ind{\sf #1}}{}
\newcommand{\pr}[1]{\ind{\it #1\/}}

\newcommand{\Section}[1]{\bigskip\section{#1}\medskip}

\newcommand{\coker}{\operatorname{coker}}
\newcommand{\im}{\operatorname{im}}
\newcommand{\id}{{\operatorname{id}}}
\newcommand{\chr}{\operatorname{char}}
\newcommand{\Hom}{\operatorname{Hom}}
\newcommand{\Tor}{\operatorname{Tor}}
\newcommand{\Ext}{\operatorname{Ext}}
\newcommand{\Cot}{\operatorname{Cot}}

\newcommand{\Nilp}{\operatorname{Nilp}}
\newcommand{\Br}{\operatorname{Br}}
\newcommand{\sub}{\subset}

\newcommand{\bwe}{\textstyle\bigwedge}
\newcommand{\bop}{\bigoplus\nolimits}
\newcommand{\bcap}{\bigcap\nolimits}

\newcommand{\lan}{\langle}
\newcommand{\ran}{\rangle}

\newcommand{\F}{{\mathbb F}}
\newcommand{\Z}{{\mathbb Z}}
\newcommand{\T}{{\mathbb T}}

\newcommand{\Q}{{\mathbb Q}}
\newcommand{\R}{{\mathbb R}}
\newcommand{\C}{{\mathbb C}}

\renewcommand{\a}{\alpha}
\renewcommand{\c}{\gamma}
\renewcommand{\o}{\omega}
\renewcommand{\O}{\Omega}
\renewcommand{\d}{\partial}
\newcommand{\e}{\varepsilon}
\newcommand{\de}{\delta}
\renewcommand{\L}{{\mathfrak L}}
\newcommand{\La}{\Lambda}
\newcommand{\D}{\Delta}
\newcommand{\G}{\Gamma}
\newcommand{\Gal}{\operatorname{Gal}}
\newcommand{\oF}{\,\overline{\!F}}
\newcommand{\KM}{K^{\mathrm{M}}}
\newcommand{\q}{{\operatorname{q}}}
\newcommand{\gr}{{\operatorname{gr}}}
\newcommand{\ab}{{\mathrm{ab}}}
\newcommand{\opp}{{\mathrm{op}}}

\newcommand{\Fradl}{{F\.[\!\!\.\sqrt[\leftroot{-2}l^\infty]{F}\.]}}
\newcommand{\Frootone}{{F\.[\!\!\.\sqrt[\leftroot{-2}l^\infty]{\.1}\.]}}
\newcommand{\Ee}{{}'\!E}
\newcommand{\dd}{{}'\!d}

\begin{document}
\vspace{0.6cm}
\title{Koszul Property and Bogomolov's conjecture}
\author{Leonid Positselski}
\address{Independent University of Moscow}
\email{posic@mccme.ru}
\maketitle

\section{Introduction}
\medskip

\subsection{}
 Let $F$ be an arbitrary field and $G_F=\Gal(\oF/F)$ be the Galois group
of its (separable) algebraic closure~$\oF$ over it.
 Two conjectures about the homological properties of the group~$G_F$
are widely known.
 First of them, the {\it Milnor--Kato conjecture}, claims that
for any prime number $l\ne\chr F$ the algebra of Galois cohomology 
with cyclotomic coefficients
$$
 H^*(G_F,\.\mu_l^{\ot *}) \.=\.
 \textstyle\bop_{n=0}^\infty \displaystyle
 H^n(G_F,\.\mu_l^{\ot n})
$$
is generated by~$H^1$ and defined by quadratic relations.
 Here, as usually, we denote by $\mu_l$ the group of $l$\+roots
of unity in~$\oF$.

 More precisely, there is a canonical homomorphism of graded algebras
called the {\it Galois symbol}, or the {\it norm residue homomorphism}
  $$
   \KM_*(F)\ot\Z/l\lrarrow H^*(G_F,\.\mu_l^{\ot *})
  $$
from the Milnor K-theory ring $\KM_*(F)$ of the field~$F$ modulo~$l$
to the cyclotomic Galois\- cohomology.
 One defines this homomorphism in degree~1 as the classical Kummer
isomorphism $F^*/F^{*l}\simeq H^1(G_F,\mu_l)$, then shows that it
can be extended to the higher degrees by multiplicativity~\cite{BT}.
 The conjecture says that for any field~$F$ and any $l\ne\chr F$
this map is an isomorphism~\cite{Mil,Kato}.

 A.~Merkurjev and A.~Suslin proved that it is an isomorphism
in degree~$n=2$; later the same authors and, independently,
M.~Rost found a proof for $l=2$ and $n=3$.
 Recently V.~Voevodsky~\cite{Voev} obtained a general proof for~$l=2$.

\subsection{}
 The second conjecture is due to F.~Bogomolov~\cite{Bog}.
 It claims that, whenever the field~$F$ contains an algebraically
closed subfield, the commutator subgroup of the Silow pro-$l$-subgroup
of~$G_F$ is a free pro-$l$-group.
 Furthermore, the commutator subgroup of the maximal quotient
pro-$l$-group of $G_F$ should be a free pro-$l$-group as well.
 The latter statement is stronger, because the former one
can be deduced from it by passing from~$F$ to the field
corresponding to the Sylow subgroup of~$G_F$.

 The condition about algebraically closed subfield cannot be dropped,
though it appears that the only essential part of it is that
the field~$F$ should contain all the roots of unity of degrees
equal to powers of~$l$.
 Indeed, in the appendix we construct simple
counterexamples showing that both the above formulations
of the conjecture fail for fields without these roots of unity.
 Nevertheless, the author thinks that there is a way to extend
Bogomolov's conjecture to all fields.

\begin{thm}{Conjecture 1}
 Let $F$ be an arbitrary field and $l$ be a prime number.
 Let $K=\Fradl$ be the field obtained by adjoining to~$F$
all roots of degrees\/~$l^N$, $\.N>0$, of all the elements of~$F$.
 Then the Sylow pro-$l$-subgroup of the absolute Galois
group~$G_K$ of the field~$K$ is a free pro-$l$-group.
\end{thm}

 It is well-known that the statement of Conjecture~1 is true
for any field of characteristic $l$~\cite{Ser}; in the sequel
we will presume that we are not in this case.
 Besides, the conjecture is true for all number fields and
their functional analogues (where it is even sufficient
to adjoin the roots of unity).
 Conjecture~1 holds for a complete discrete valuation field
whenever it holds for the residue field.
 This is the most important evidence supporting Conjecture~1
that is currently known to the author.

\subsection{}
 Conjecture~1 is very strong and likely hard to approach.
 The following weaker version is closer to original Bogomolov's
conjecture.

\begin{thm}{Conjecture 2}
 Let $F$ be  a field containing a primitive root of unity of degree~$l$
if\/ $l$ is odd, and a field containing a square root of $-1$ if\/ $l=2$.
 Then the maximal quotient pro-$l$-group of the absolute Galois
group\/~$G_K$ of the field\/ $K=\Fradl$ is a free pro-$l$-group.
\end{thm}

 Conjecture~2 follows from Conjecture~1, because the maximal
quotient pro-$l$-group $G^{(l)}$ of a pro-finite group~$G$ is free
whenever the Sylow pro-$l$-subgroup $S_l\.G$ of~$G$ is free.
 Indeed, both maps in the sequence $H^2(G^{(l)},\.\Z/l)\rarrow
H^2(G,\.\Z/l)\rarrow H^2(S_l\.G,\.\Z/l)$ are easily seen to be
injective for any pro-finite group~$G$.

\subsection{}
 The relation between the Milnor--Kato and freeness conjectures
is the main topic of Bogomolov's paper~\cite{Bog}.
 It seems that no one of the two conjectures implies directly
the other one; rather, they are somewhat complementary.

\subsection{}
 The approach to the Milnor--Kato and freeness conjectures which we
develop here was started in the paper~\cite{PV} by Alexander Vishik
and me.
 There it was proven that the statement of the Milnor--Kato conjecture
for a field~$F$ whose Galois group $G_F$ is a pro-$l$-group would follow
from its low-degree part if we knew the Milnor K\+theory ring
$\KM_*(F)\ot\Z/l$ to be a {\it Koszul algebra}.
 The following theorem, which we prove in section~5 (as the first
statement of Corollary~2) slightly generalizes the main result
of~\cite{PV}.

\begin{thm}{Theorem 1}
 Let $l$ be a prime number and $F$ be a field containing a primitive
$l$\+root of unity.
 Suppose that
 \begin{itemize}
   \item[(1)] the map $\KM_n(F)\ot\Z/l\rarrow H^n(G_F,\.\mu_l^{\ot n})$
          is an isomorphism for~$n=2$ and a monomorphism for~$n=3$;
   \item[(2)] the algebra $\KM_*(F)\ot\Z/l$ is Koszul.
 \end{itemize}
 Let $G_F^{(l)}$ be the maximal quotient pro-$l$-group of
the group~$G_F$.
 Then the natural homomorphism
$\KM_*(F)\ot\Z/l\rarrow H^*(G_F^{(l)},\.\mu_l^{\ot *})$
is an isomorphism.%
\end{thm}

 Of course, the first part of the condition~(1) is the theorem of
Merkurjev and Suslin; the second part of~(1) is long known to be true,
at least, for $l=2$.
 It is not difficult to show, using the transfer operations on
the Milnor K-theory, that it suffices to verify the Milnor--Kato
conjecture in the pro-$l$-group case~\cite{Bog,SV}.
 Thus our Koszulity hypothesis essentially implies the Milnor--Kato.

\subsection{}
 The main goal of this paper is to show that a certain much stronger
hypothesis about the Milnor ring would imply Conjecture~2
(Bogomolov's conjecture) as well.

 Under the conditions of Conjecture~2, the Milnor K-theory algebra
modulo~$l$ is a quotient algebra of the exterior algebra of
the vector space $F^*/F^{*l}$ over the prime field $\Z/l$, so we
can consider the natural homomorphism
  $$
    \Lambda_l^*(F) =
    \bwe_{\Z/l}^*(F^*/F^{*l}) \lrarrow \KM_*(F)\ot\Z/l.
  $$
 Notice that the exterior algebra $\Lambda_l^*(F)$ is isomorphic to
the cohomology algebra of the Galois group $\Gal(\Fradl/F)$ with
coefficients $\Z/l$ (see section~9).

 The stronger Koszulity hypothesis states that the kernel $J_l^*(F)$
of this map should be a {\it Koszul module\/} over the exterior
algebra (see section~3 for the definition).
 Our main result can be formulated as follows.

\begin{thm}{Theorem 2}
 Let $l$ be a prime number and $F$ be a field containing the\/ $l$\+roots
of unity if\/ $l$~is odd, and containing the $4$\+roots of unity if\/
$l=2$.
 Suppose that
 \begin{itemize}
   \item[(1)] the condition~(1) of Theorem~1 holds;
   \item[(2)] the ideal $J_l^*(F)$ generated by the Steinberg relations
       in the exterior algebra $\Lambda_l^*(F)$ is a Koszul module
       (in the grading shifted by~$2$) over the algebra $\Lambda_l^*(F)$.
 \end{itemize}
 Then the maximal quotient pro-$l$-group of the absolute Galois
group $G_K$ of the field $K=\Fradl$ is a free pro-$l$-group.
\end{thm}

 It is shown in section~6 that the condition~(2) of Theorem~1 follows
from the condition~(2) of Theorem~2.
 Therefore, in the assumptions of Theorem~2 both the Milnor--Kato
and Bogomolov's conjectures are true.

\subsection{}
 Two proofs of Theorem~2 are given in this paper.
 The first one is rather computational; it is based on the
Serre--Hochschild spectral sequence for pro-$l$-groups and
a result extending the techniques of~\cite{PV} to cohomology
of conilpotent comodules.
 The latter is proven in section~4.
 The second argument uses the results of~\cite{PV} in order to pass
to the associated graded coalgebras and then applies directly
a general principle claiming that Koszulity-type conditions
on morphisms of Koszul algebras correspond to freeness-type
conditions under duality.
 This is explained in section~7.

\subsection{}
 Finally, let us point to a simple application of our results
to a different area of motivic theory---that of mixed Tate motives
with rational coefficients over a field of finite characteristic.
 For such a field $F$, the standard conjectures~\cite{Beil}
imply that the Milnor K-theory algebra with rational coefficients
$\KM_*(F)\ot\Q$ should coincide with the algebraic K-theory
$K_*(F)\ot\Q$ and be a Koszul algebra.
 Furthermore, in this case the {\it motivic pro-Lie algebra\/}
$\L_*(F)=\bop_{i=1}^\infty \L_i(F)$ describing the abelian tensor
category of mixed Tate motives with rational coefficients
over a field~$F$ is isomorphic to the quadratic pro-Lie algebra
dual to $\KM_*(F)\ot\Q$.
 (The duality here connects the quadratic pro-Lie algebras with
the skew-commutative quadratic algebras.)

 A conjecture of A.~Goncharov~\cite{Gon} claims that for any
field~$F$ the subalgebra $\L_{\ge2}(F)=\bop_{i=2}^\infty \L_i(F)$ of
the motivic pro-Lie algebra~$\L_*(F)$ is a free graded pro-Lie algebra.
 It follows from our result that this conjecture together with
the Koszulity conjecture above is equivalent to the following
hypothesis about the Milnor K-theory algebra.
 Let $J_\Q^*(F)$ be the kernel of the homomorphism from the exterior
algebra $\bwe^*_\Q(F^*{\ot}\.\Q)$ to $\KM_*(F)\ot\Q$; then the ideal
$J_\Q^*(F)$ should be a Koszul module over the exterior algebra.
 Note that in this setting, unlike in the Bogomolov case, the analogue
of the statement converse to Theorem~2 also holds: the algebra
$\KM_*(F)\ot\Q$ is Koszul and the pro-Lie algebra $\L_{\ge2}(F)$
is free if and only if the ideal $J_\Q^*(F)$ is a Koszul module.
 The basic reason is that $\L_*(F)$ is a graded pro-Lie algebra,
while in the Bogomolov case we deal with an ungraded Galois group.
 In this sense, the hypothesis of Koszulity of the Steinberg
ideal $J_l^*(F)$ can be thought of as a stronger graded version
of Conjecture~2.
 This is explained in section~8.

\subsection{}
 Coalgebras, comodules, and their cohomology are considered
in section~2.
 Definitions and basic properties of Koszul (co)algebras and
(co)modules are presented in section~3.
 The cohomology of conilpotent comodules are studied in section~4,
the cohomology of nonconilpotent coalgebras in section~5.
 We discuss Koszul properties related to morphisms of graded
algebras in section~6.
 The duality between Koszulity and freeness conditions is
explained in section~7.
 Koszulity and freeness for Lie coalgebras are considered
in section~8.
 We prove the main theorem in section~9.
 Evidence related to Conjecture~1 is presented in the appendix.

\bigskip
 The author wishes to express his gratitude to V.~Voevodsky for posing
the problem and to A.~Beilinson, F.~Bogomolov, A.~Goncharov, F.~Pop,
and A.~Yakovlev for very helpful discussions.
 I owe to my advisor B.~Mazur the idea of the counterexample in
the appendix, from which the above formulation of Conjectures~1\+-2
eventually grew.
 Most of the mathematical content of this paper was invented when
the author was a graduate student at Harvard University (the paper
is an enhanced version of my Ph.~D. thesis).
 The rest was done when I was working at the Independent University
of Moscow and visiting the Max-Planck-Institut f\"ur Mathematik in Bonn.
 The author wishes to thank all the mentioned institutions.

\Section{Homology of Algebras and Cohomology of Coalgebras}

\subsection{The cohomology of comodules}
 Recall that a {\it coalgebra\/} $C$ over a field~$\k$ is a vector space
equipped with a comultiplication map $\D\:C\rarrow C\ot C$ and
a counit map $\e\:C\rarrow\k$ satisfying the conventional coassociativity
and counit axioms.
 A {\it coaugmented coalgebra\/} is a coalgebra $C$ endowed with a
coalgebra homomorphism $\c\:\k\rarrow C$.
 The cohomology algebra $H^*(C)$ of a coaugmented coalgebra $C$
is the opposite multiplication algebra to the Ext-algebra
$\Ext^*_C(\k,\k)$, where~$\k$ is endowed with the left $C$\+comodule
structure by means of~$\c$.

 A left {\it comodule\/} $P$ over a coalgebra $C$ is a vector space
over~$\k$ together with a coaction map $\D'\:P\rarrow C\ot P$
satisfying the usual coaction axiom.
 The cohomology module of a comodule $P$ over a coaugmented coalgebra~$C$
is the left $H^*(C)$\+module $H^*(C,P)=\Ext^*_C(\k,P)$.
 This cohomology can be calculated in terms of the following
explicit cobar-resolution
 $$
  P\lrarrow C\ot P\lrarrow C\ot C^+\ot P\lrarrow C\ot C^+\ot C^+\ot P
  \lrarrow \dsm,
 $$
where $C^+=\coker\c$ and the the cobar-differential is given by
the well-known formula $d(c_0\ot\ds\ot c_n\ot p)=
\sum_{i=0}^n (-1)^{i-1} c_0\ot\ds\ot\D(c_i)\ot\ds\ot c_n\ot p
+(-1)^n c_0\ot\ds\ot c_n\ot \D'(p)$.
 It is easy to check that this is an injective comodule resolution
of~$P$.
 Applying the functor $\Hom_C(\k,{-})$ to it, we get
 $$
   H^*(C,P)=H^*(P\rarrow C^+\ot P\rarrow C^+\ot C^+\ot P\rarrow\dsb)
 $$
and in particular
 $$
   H^*(C)=H^*(\.\k\rarrow C^+\rarrow C^+\ot C^+\rarrow\dsb).
 $$
 The multiplication on $H^*(C)$ and the action of $H^*(C)$ on
$H^*(C,P)$ are induced by the obvious multiplication and action
 \begin{align*}
   (c_1\ot\ds\ot c_i)\cd(c_{i+1}\ot\ds\ot c_{i+j})&=
   c_1\ot\ds\ot c_i\ot c_{i+1}\ot\ds\ot c_{i+j}         \\
   (c_1\ot\ds\ot c_i)\cd(c_{i+1}\ot\ds\ot c_{i+j}\ot p)&=
   c_1\ot\ds\ot c_i\ot c_{i+1}\ot\ds\ot c_{i+j}\ot p
 \end{align*}
on the level of cobar-complexes.

\subsection{The homology of modules}
 The homology coalgebra $H_*(A)$ of an augmented algebra~$A$
is defined as $H_*(A)=\Tor^A_*(\k,\k)$, where the left and
the right $A$\+module structures on~$\k$ are defined by means
of the augmentation morphism $\a\:A\rarrow\k$.
 Let $M$ be a left module over an augmented algebra~$A$.
 The homology comodule of~$M$ is the left $H_*(A)$\+comodule
$H_*(A,M)=\Tor^A_*(\k,M)$.
 This homology can be calculated by means of the bar-resolution
  $$
   M\llarrow A\ot M\llarrow A\ot A_+\ot M
  \llarrow A\ot A_+\ot A_+\ot M\llarrow\dsb
  $$
with $A_+=\ker\a$ and the standard bar-differential.
 Therefore we have
 $$
  H_*(A,M)=H_*(M\larrow A_+\ot M\larrow A_+\ot A_+\ot M\larrow\dsb)
 $$
and in particular
 $$
  H_*(A)=H_*(\.\k\larrow A_+\larrow A_+\ot A_+\larrow\dsb).
 $$
 The comultiplication on $H_*(A)$ and the coaction of $H_*(A)$
on $H_*(A,M)$ are induced by the obvious comultiplication and coaction
 \begin{align*}
   \D(a_1\ot\ds\ot a_n)&=\sum_{i=0}^n
   (a_1\ot\ds\ot a_i)\ot(a_{i+1}\ot\ds\ot a_n)  \\
   \D'(a_1\ot\ds\ot a_n\ot m)&=\sum_{i=0}^n
   (a_1\ot\ds\ot a_i)\ot(a_{i+1}\ot\ds\ot a_n\ot m)
 \end{align*}
on the level of bar-complexes.

\subsection{Graded modules and comodules}
 By a {\it graded coalgebra\/} we mean a nonnegatively graded
vector space $C=\bop_{n=0}^\infty C_n$ equipped with a coalgebra
structure which respects the grading, i.~e.,
$\D(C_n)\sub\sum_{i+j=n}C_i\ot C_j$ and $\e(C_{>0})=0$.
 Moreover, we will always assume the map $\e:C_0\rarrow\k$
to be an isomorphism.
 Then a graded coalgebra structure induces a coaugmented coalgebra
structure in the obvious way.
 A {\it graded comodule\/} over $C$ is a graded vector space
$P=\bop_{i=0}^\infty P_i$ equipped with a $C$\+comodule structure
which respects the grading,
i.~e., $\D'(P_n)\sub\sum_{i+j=n}C_i\ot P_j$.
 We always consider the nonnegatively graded comodules only,
that is, $P_i=0$ for $i<0$.

 We leave it to the reader to formulate the dual definitions
of graded algebras and modules; as above, we assume that
$A_i=0$ for $i<0$, that $A_0$ is a one-dimensional vector space
generated by the unit, and that $M_i=0$ for $i<0$.

 One can immediately see from the bar-complexes above that
the homology and cohomology of graded objects are endowed with
a natural second grading: for a graded algebra $A$ and
a graded $A$\+module $M$ we have
  $$
   H_*(A)=\bop_{i\le j}H_{i,j}(A)
   \quad \text{and} \quad
   H_*(A,M)=\bop_{i\le j}H_{i,j}(A,M);
  $$
and analogously, for a graded coalgebra $C$ and
a graded $C$\+comodule $P$
  $$
   H^*(C)=\bop_{i\le j}H^{i,j}(C)
   \quad \text{and} \quad
   H^*(C,P)=\bop_{i\le j}H^{i,j}(C,P).
  $$

\subsection{Opposite multiplications and comultiplications}
 For algebras $A$ and~$A^\opp$ with opposite multiplication
there is a canonical isomorphism of graded vector spaces
$H_*(A)\simeq H_*(A^\opp)$ given by the tautological isomorphism
$\Tor^A_*(\k,\k)\simeq\Tor^{A^\opp}_*(\k,\k)$.
 It is easy to see from the bar-complexes above that this
isomorphism makes $H_*(A)$ and $H_*(A^\opp)$ coalgebras
with opposite comultiplication.

 On the other hand, for coalgebras $C$ and~$C^\opp$ with opposite
comultiplication there is no self-obvious isomorphism between
$\Ext_C^*(\k,\k)$ and $\Ext_{C^\opp}^*(\k,\k)$.
 Nevertheless, one can use the cobar-complexes to construct
an isomorphism of graded vector spaces $H^*(C)\simeq H^*(C^\opp)$
which will make them algebras with opposite multiplication.

 To define this isomorphism without mentioning explicit resolutions,
one has to introduce the {\it cotensor product\/} and {\it cotorsion\/}
functors, as explained in the next subsection.
 This will complete the symmetry between the algebras and coalgebras.

\subsection{The cotorsion functors and cohomology of comodules}
 Let $P$ and~$Q$ be a right and a left comodule over a coalgebra~$C$.
 The {\it cotensor product\/} vector space $P\oa_C Q$ is defined as
the kernel of the map $\D'\ot\id-\id\ot\D'\:P\ot Q\rarrow P\ot C\ot Q$.
 It is not difficult to check that this operation satisfies
the associativity and unit identities analogous to those for
the tensor product of modules over algebras:
$(P\oa_C Q)\oa_D R=P\oa_C(Q\oa_D R)$ and $C\oa_C P=P$.

 The cotensor product is a left-exact functor on both arguments;
according to the above identities, the functors of cotensor
product with cofree comodules are exact.
 As usual, it follows that the right derived functors of the functor
of cotensor product along the first and the second arguments coincide;
we will call them the {\it cotorsion functors\/}
and denote by $\Cot_C^i(P,Q)$.
 Using the cobar-resolution from subsection~1.1, one can compute
those derived functors explicitly:
 $$
  \Cot_C^*(P,Q)=
   H^*(P\ot Q\rarrow P\ot C\ot Q\rarrow P\ot C\ot C\ot Q\rarrow\dsb),
 $$
where $C$ can be replaced with~$C_+=\ker\e$.

 For any left $C$\+comodule $P$ and finite-dimensional left
$C$\+comodule $Q$ there are natural isomorphisms
$\Ext_C^i(Q,P)\simeq\Cot_C^i(Q^*,P)$, just as for any
left $A$\+module $M$ and right $A$\+module $N$ one has
$\Ext_A^i(M,N^*)\simeq\Tor^A_i(N,M)^*$.
 Therefore, we may alternatively define the cohomology of
a coaugmented coalgebra as $H^*(C)=\Cot_C^*(\k,\k)$ and
the cohomology of a comodule as $H^*(C,P)=\Cot_C^*(\k,P)$.

\begin{rem}{Remark:}
 All the definitions of this section can be given for the (co)algebra
and (co)module structures on (graded) objects of an arbitrary (not
necessarily semisimple; associative, but not necessarily commmutative)
abelian tensor category (with an exact functor of tensor product).
 The same holds for all the definitions and results of section~3,
6 and~7; all of them are duality-symmetric, in particular.
 For the results of section~8 mentioning commutative and Lie
structures, the tensor category has to be commutative, of course.
 All of this does not apply to the results of sections~4 and~5, however,
which depend essentially on the standard properties of the functor
of direct limit and therefore are not duality-symmetric at all.
 But we prefer to deal with vector spaces in this paper.
\end{rem}

\Section{Koszul Property for Modules and Comodules}

 Different accounts of the basic properties of {\it quadratic\/}
and {\it Koszul algebras\/} can be found
in the papers~\cite{Pr,BF,BGSoe,BGS}.
 The notion of a {\it Koszul module\/} and related results first
appeared in~\cite{BGSoe}. 
 In the infinite-dimensional case one has to define
the {\it quadratic duality\/} as a correspondence
between algebras and coalgebras~\cite{PV}.
 The purpose of this section is to present the analogous results
for modules and comodules.

\subsection{Quadratic modules and comodules}
 A graded algebra $A$ is called {\it one-generated\/} if it is
multiplicatively generated by $A_1$.
 A graded algebra $A$ is called {\it quadratic\/} if it is isomorphic
to the quotient algebra $\{V,R\}=\T(V)/(R)$ of the tensor algebra
$\T(V)=\bop_n V^{\ot n}$ of the vector space $V=A_1$
by the ideal generated by a subspace $R\sub V^{\ot2}$.
 For any graded algebra $A$ there exists a unique quadratic algebra
$\q A$ together with a morphism of graded algebras $r_A\:\q A\rarrow A$
which is an isomorphism on $A_1$ and a monomorphism on $A_2$.
 To construct the algebra $\q A$, one can take $A_1$ for~$V$ and
the kernel of the multiplication map $A_1\ot A_1\rarrow A_2$ for~$R$.

 A graded module $M$ over a one-generated algebra $A$ is called
{\it one-generated\/} if it is generated by $M_0$, that is
the map $A\ot M_0\rarrow M$ is surjective.
 A graded module $M$ over a quadratic algebra $A=\{V,R\}$ is called
{\it quadratic\/} if it is isomorphic to the quotient module
$\{U,S\}_A=A\ot U/AS$ of the free $A$\+module $A\ot U$ generated by
the vector space $U=M_0$ by the submodule generated by a vector subspace
$S\sub V\ot U$.
 For any graded module $M$ over a graded algebra $A$ there exists
a unique quadratic module $\q_AM$ over the quadratic algebra $\q A$
together with a morphism of graded $\q A$\+modules
$r_{A,M}\:\q_AM\rarrow M$ which is an isomorphism on $M_0$
and a monomorphism on $M_1$.
 To construct the module $\q_AM$, one takes $M_0$ for~$U$ and
the kernel of the map $A_1\ot M_0\rarrow M_1$ for~$S$.

 A graded coalgebra $C$ is called {\it one-cogenerated\/}
if the iterated comultiplication maps
$\D^{(n)}\:C_n\rarrow C_1^{\ot n}$ are injective, or equivalently,
all the maps $\D\:C_{i+j}\rarrow C_i\ot C_j$ for $i$,~$j\ge0$
are injective.
 A graded coalgebra is called {\it quadratic\/} if it is isomorphic
to the subcoalgebra of the tensor coalgebra $\bop_n V^{\ot n}$
of the form
 $$
  \lan V,R\ran=\bop_{n=0}^\infty\,
  \bcap_{i=1}^{n-1}\. V^{i-1}\ot R\ot V^{n-i-1}
 $$
for the vector space $V=C_1$ and a certain subspace $R\sub V^{\ot2}$.
 For any graded coalgebra $C$ there exists a unique quadratic coalgebra
$\q C$ together with a morphism of graded coalgebras $r_C\:C\rarrow \q C$
which is an isomorphism on $C_1$ and an epimorphism on $C_2$.
 Indeed, it suffices to take $C_1$ for~$V$ and the image
of the comultiplication map $C_2\rarrow C_1\ot C_1$ for~$R$.

 A graded comodule $P$ over a one-cogenerated coalgebra $C$
is called {\it one-cogenerated\/} (cogenerated by $P_0$)
if all the coaction maps $\D'\:P_{i+j}\rarrow C_i\ot P_j$
for $i$,~$j\ge0$ are injective, or equivalently, the maps
$\D'\:P_i\rarrow C_i\ot P_0$ are injective.
 A graded comodule $P$ over a quadratic coalgebra $C=\lan V,R\ran$
is called {\it quadratic\/} if it is isomorphic to the subcomodule
of the cofree $C$\+comodule $C\ot U$ of the form
 $$
   \lan U,S\ran_C=\bop_{n=0}^\infty\. C_n\ot U\cap C_{n-1}\ot S
 $$
for the vector space~$U=P_0$ and a certain subspace $S\sub V\ot U$.
 For any graded comodule~$P$ over a graded coalgebra~$C$ there exists
a unique quadratic comodule $\q_CP$ over the quadratic coalgebra
$\q C$ together with a morphism of graded $\q C$\+comodules
$r_{C,P}\:P\rarrow\q_CP$ which is an isomorphism on $P_0$
and an epimorphism on $P_1$.
 Namely, one can take $P_0$ for~$U$ and the image of the coaction map
$P_1\rarrow C_1\ot P_0$ for~$S$.

\subsection{Quadratic duality}
 The quadratic algebra $A=\{V,R\}$ and the quadratic coalgebra
$C=\lan V,R\ran$ are called {\it dual\/} to each other; we denote
this as $A=C^!$ and $C=A^?$.
 This rule defines an equivalence between the category of
quadratic algebras and the category of quadratic coalgebras.
 Furthermore, the quadratic $A$\+module $M=\{U,S\}_A$ and the quadratic
$C$\+comodule $P=\lan U,S\ran_C$ are called {\it dual\/} to each other,
too; the notation: $M=P_C^!$ and $P=M_A^?$. 
 This defines an equivalence between the category of quadratic
$A$\+modules and the category of quadratic $C$\+comodules.

\begin{thm}{Proposition 1}
 A graded coalgebra $C$ is one-cogenerated if and only if
$H^{1,j}(C)=0$ for all $j>1$.
 A one-cogenerated coalgebra $C$ is quadratic if and only if
$H^{2,j}(C)=0$ for all $j>2$.
 More precisely, for any graded coalgebra~$C$ the morphism
$r_C\:C\rarrow\q C$ is an isomorphism in degree\/~$\.\le n$
if and only if $H^{i,j}(C)=0$ for all\/ $i<j\le n$ and\/
$i=1$,~$2$.
 The analogous statements are true for graded algebras
and their homology.%
\end{thm}

\pr{Proof}: See~\cite{PV}.  \qed

\begin{thm}{Proposition 2}
 A graded comodule $P$ over a one-cogenerated coalgebra $C$
is one-cogenerated if and only if $H^{0,j}(C,P)=0$ for all $j>0$.
 A one-cogenerated comodule $P$ over a quadratic coalgebra $C$
is quadratic iff $H^{1,j}(C,P)=0$ for all $j>1$.
 More precisely, for any graded comodule~$P$ over a quadratic
coalgebra~$C$ the morphism $r_{C,P}\:P\rarrow\q_C P$
is an isomorphism in degree\/~$\.\le n$ if and only if
$H^{i,j}(C,P)=0$ for all\/ $i<j\le n$ and\/ $i=0$,~$1$.
 The analogous statements are true for graded modules.
\end{thm}

\pr{Proof}:
 The argument is based on the explicit form of the cobar-complex
computing the cohomology in question.
 The space $H^0(C,P)$ is isomorphic to the kernel of the
map $\D'\:P\rarrow C_+\ot P$; if the coalgebra $C$
is one-cogenerated, it follows easily that the morphism
$r_{C,P}\:P\rarrow \q_CP\sub C\ot P_0$ is injective
in degree~$\le n$ if and only if $H^{0,j}(C,P)=0$ for all $0<j\le n$.
 Now assume that the coalgebra~$C$ is quadratic and the map~$r_{C,P}$
is injective in degree $\.<\nbk n$.
 Let $z\in C^+\ot P$ be a homogeneous cocycle of degree~$n$;
then $z=\sum_{s\ge1,t\ge0}^{s+t=n}z_{st}$ with $z_{st}\in C_s\ot P_t$.
 The cocycle condition means that the images of $(\D\ot\id)(z_{u+v,w})$
and $(\id\ot\D')(z_{u,v+w})$ in $C_u\ot C_v\ot P_w$ coincide
for any $u,v\ge1$, \ $w\ge0$, \ $u+v+w=n$.
 Since the comultiplication and coaction maps are injective,
the latter equivalently means that the images of $z_{st}$ in
$C_1^{\ot n}\ot P_0$ coinside for all $s$,~$t$.
 So we get an element of $C_1^{\ot n}\ot P_0$;
it is easy to see that it represents an element of $\q_CP$.
 Besides, this element belongs to the image of $r_{C,P}$ iff
the cocycle~$z$ is a coboundary.
 Finally, if the map $r_{C,P}$ is an isomorphism in degree~$\.<\nbk n$,
then any element of~$\q_CP$ of degree~$n$ corresponds to a cocycle~$z$
in this way.
\qed

\begin{thm}{Proposition 3}
 For any graded coalgebra $C$, the diagonal subalgebra\/
$\bop_i H^{i,i}(C)$ of the cohomology algebra $H^*(C)$ is
a quadratic algebra isomorphic to $(\q C)^!$.
 For any graded $C$\+comodule $P$, the diagonal part\/
$\bop_i H^{i,i}(C,P)$ of the cohomology module $H^*(C,P)$
is a quadratic module over the diagonal subalgebra of $H^*(C)$
isomorphic to the $(\q C)^!$\+module $(\q_CP)_{\q C}^!$.
 Analogously, for any graded algebra $A$ the diagonal\-
homology coalgebra\/ $\bop_i H_{i,i}(A)$ is isomorphic
to the quadratic coalgebra~$(\q A)^?$.
 For any graded $A$\+module $M$, the diagonal homology\/
$\bop_i H_{i,i}(A,M)$ is a quadratic comodule over
the coalgebra\/ $\bop_i H_{i,i}(A)$ isomorphic to
the comodule $(\q_AM)_{\q A}^?$ over $(\q A)^?$.
\end{thm}

\pr{Proof}: Straightforward computation with the bar
and cobar-complexes. \qed

\subsection{Koszul modules and comodules}
 The notion of a {\it Koszul algebra\/} and the construction
of the {\it Koszul complex\/} were invented by S.~Priddy
in~\cite{Pr}.
 The result relating Koszulity with distributivity of collections 
of vector subspaces is due to J.~Backelin~\cite{Bac}.
 Detailed expositions can be found in the papers~\cite{BGSoe,BGS};
see also~\cite{BF}.

 A graded coalgebra $C$ is called {\it Koszul\/} if one has
$H^{i,j}(C)=0$ for all $i\ne j$.
 A graded comodule $P$ over a Koszul coalgebra $C$ is called
{\it Koszul\/} if $H^{i,j}(C,P)=0$ for $i\ne j$.
 Analogously, a graded algebra $A$ is called {\it Koszul\/}
if $H_{i,j}(A)=0$ unless $i=j$.
 A graded module $M$ over a Koszul algebra $A$ is called {\it Koszul\/}
if $H^{i,j}(A,M)=0$ unless $i=j$.
 Note that all Koszul algebras and coalgebras, as well as modules
and comodules, are quadratic.
 This follows from Propositions~1 and~2.

 It is clear from subsection~1.4 that algebras with opposite
multiplication, as well as coalgebras with opposite comultiplication,
are quadratic or Koszul simultaneously.
 A graded right module over a Koszul algebra~$A$ is called {\it Koszul\/}
if it is a Koszul module over $A^\opp$; analogously for right comodules.

 Let $A$ be a quadratic algebra and $A^?$ be the dual coalgebra;
suppose we are given a right $A$\+module $M$ and a left
$A^?$\+comodule $P$.
 The {\it Koszul complex\/} $K(M,P)$ is defined as the vector
space $M\ot P$ with the differential given by the formula
 $$
  \d=(m'\ot\id)(\id\ot \D')\:M\ot P\lrarrow M\ot A_1\ot P\lrarrow M\ot P,
 $$
where $m'\: M\ot A_1\rarrow M$ is the action morphism.
 It follows immediately from the definition of~$A^?$ that one has
$\d^2=0$.
 For a left $A$\+module $M$ and a right $A^?$\+comodule $P$, the
{\it Koszul complex\/} $K(P,M)$ is defined in the analogous way.

\begin{thm}{Proposition 4}
 A quadratic algebra~$A$ and the dual quadratic coalgebra~$A^?$
are Koszul simultaneously.
 Moreover, they are Koszul if and only if the Koszul complex
$K(A,A^?)$ is exact at all its terms $A_p\ot A^?_q$ with
$p+q>0$.
 A quadratic module~$M$ over a Koszul algebra~$A$ and the dual
quadratic comodule~$M^?_A$ over the dual Koszul coalgebra~$A^?$
are Koszul simultaneously, too.
 Furthermore, they are both Koszul if and only if the Koszul
complex $K(A,M^?_A)$ is a resolution of the $A$\+module $M$ and
if and only if the Koszul complex $K(A^?,M)$ is a resolution
of the $A^?$\+comodule $M^?_A$.
\end{thm}

\pr{Proof}:
 Let us introduce a homological grading on the complex $K(A,M^?_A)$
by the rule $K_i(A,M^?_A)=A\ot M^?_{A,\,i}$.
 There is a natural morphism of complexes of graded $A$\+modules
from $K(A,M^?_A)$ to the bar-resolution
$A\ot M\larrow A\ot A_+\ot M\larrow\dsb$ of the $A$\+module $M$
defined by the formula $\rho(a\ot p_i)=a\ot\D^{\1\,(i)}(p_i)$,
where $\D^{\1\,(i)}\: M^?_{A,\,i}\rarrow
A_1^{\ot i}\ot M_0\sub A_+^{\ot i}\ot M$ is the natural embedding
given by the iterated coaction map.
 The composition of~$\rho$ with the standard morphism from
the bar-resolution to the module~$M$ itself provides a morphism
$K(A,M^?_A)\rarrow M$.

 Note that both the Koszul complex and the bar-resolution are complexes
of free graded $A$\+modules.
 It is not difficult to see (a version of Nakayama's lemma) that
a morphism of complexes of free graded modules over an arbitrary
graded algebra~$A$ is a quasi-isomorphism if and only if it becomes
a quasi-isomorphism after tensoring with~$\k$ over~$A$, provided that
both complexes are concentrated in the positive homological degree.
 Tensoring with~$\k$ turns the bar-resolution to the bar-complex
computing $H_*(A,M)$ and the complex $K(A,M)$ to the complex with
the terms $M^?_{A,\,i}$ and zero differentials; now it follows from
Proposition~3 that the map $\k\ot_A\rho$ is a quasi-isomorphism
iff $H_{i,j}(A,M)=0$ for $i\ne j$.

 Applying this result to the trivial $A$\+module $M=\k$, we conclude
that a quadratic algebra~$A$ is Koszul if and only if the Koszul
complex $K(A,A^?)$ is a resolution of~$\k$.
 Besides, we proved that a quadratic module~$M$ over a Koszul
algebra~$A$ is Koszul iff the complex $K(A,M^?_A)$ is a resolution
of the module $M$.
 The same arguments as above applied to a morphism of complexes
of cofree graded comodules show that the comodule~$M^?_A$ over~$A^?$
is Koszul iff the complex $K(A^?,M)$ is a resolution of~$M^?_A$.

 On the other hand, since the algebra $A^\opp$ is Koszul, the complex
$K(A^?,A)$ is a free graded resolution of the trivial right $A$\+module.
 Computing the spaces $\Tor^A_*(\k,M)$ in terms of this resolution
of~$\k$, we get $K(A^?,A)\ot_AM=K(A^?,M)$ and $H_*(A,M)=H_*K(A^?,M)$.
 So it follows from Proposition~3 that the map $K(A^?,M)\rarrow M^?_A$
is an isomorphism on the homology iff the $A$\+module $M$ is Koszul.
 Now we see that the $A$\+module $M$ and the $A^?$\+comodule $M^?_A$
are Koszul simultaneously.
\qed

\begin{rem}{Remark:}
 The techniques of the above argument allow to figure out quite
precisely the relations between the particular pieces of the
Koszulity condition for the dual quadratic algebras and modules.
 For example, given two integers $a$,~$b\ge2$ one has
$H_{i,j}(A)=0$ for all $i\le a+1$ and $0<j-i\le b-1$
if and only if $H^{i,j}(A^?)=0$ for all $i\le b+1$
and $0<j-i\le a-1$.
 Moreover, there are natural morphisms of vector spaces
$H^{b+1,a+b}(A^?)\rarrow H_{a+1,a+b}(A)$, which annihilate
the nontrivial (co)mul\-ti\-pli\-ca\-tions on both sides.
 There are analogous results for dual quadratic modules.
\end{rem}

\subsection{Koszulity and distributivity}
 A collection of subspaces $X_1$, \ds, $X_{n-1}$ in a vector space
$W$ is called {\it distributive\/} if there exists a (finite)
direct decomposition $W=\bop_{\o\in\O}W_\o$ such that each
subspace $X_k$ is the sum of a set of subspaces $W_\o$.
 Equivalently, a collection of subspaces~$X_k$ is distributive if
the distributivity identity $(X+Y)\cap Z=X\cap Z+Y\cap Z$
is satisfied for any triple of subspaces $X$, $Y$, $Z$ which can be
obtained from the subspaces $X_k$ using the operations of sum and
intersection.

\begin{thm}{Proposition 5}
 The quadratic algebra $A=\{V,R\}$ and the quadratic coalgebra
$C=\lan V,R\ran$ are Koszul if and only if the collection of subspaces
 $$
  V^{\ot k-1}\ot R\ot V^{n-k-1}\;\sub\; V^{\ot n}, \qquad k=1,\ds,n-1
 $$
is distributive for all\/ $n\ge 4$.
 The quadratic module $\{U,S\}_A$ over a Koszul algebra~$A=\{V,R\}$
and the quadratic comodule $\lan U,S\ran_C$ over the dual Koszul
coalgebra~$C=\lan V,R\ran$ are Koszul if and only if the collection
of subspaces
 $$
  V^{\ot k-1}\ot R\ot V^{n-k-1}\ot U, \quad k=1,\ds,n-1; \qquad
  V^{\ot n-1}\ot S
 $$
in the vector space $V^{\ot n}\ot U$ is distributive
for all\/ $n\ge 3$.
\end{thm}

\pr{Proof}:
 The first statement is proven in~\cite{PV}; the proof of the
second one is analogous and based on the same lemma. \qed

\begin{rem}{Remark:}
 Many natural Koszulity-type homological conditions on data
involving Koszul algebras and modules (including most of the
properties that will be considered in section~5) can be simply
rewritten in terms of distributivity of various collections
of subspaces.
 It is not at all clear why this so happens.
 The related results often can be proven by both homological
and lattice-theoretical methods, though the former are generally
more powerful.
 Here are two of the various examples.

 Let $A=\{V,R\}$ be a Koszul algebra and $E\sub V$ be a vector
subspace.
 Then the subalgebra generated by~$E$ in the algebra~$A$ is Koszul
if and only if the collection of subspaces $E^{\ot n}$ and
$V^{\ot k-1}\ot R\ot V^{n-k-1}$ in the vector space $V^{\ot n}$
is distributive for all $n\ge 3$.
 Note that in general this subalgebra does not have to be quadratic.

 Furthermore, let $N=W\ot A/TA$ and $M=A\ot U/AS$, where
$T\sub W\ot V$ and $S\sub V\ot U$, be a left and a right
Koszul $A$\+modules.
 Then one has $\Tor^A_{i,n}(N,M)=0$ for all $i\ne 0$,~$n$
if and only if the collection of subspaces 
$W\ot V^{\ot k-1}\ot R\ot V^{n-k-1}\ot U$ together with
$W\ot V^{\ot n-1}\ot S$ and $T\ot V^{\ot n-1}\ot U$
is distributive in $W\ot V^{\ot n}\ot U$.
 The latter observation was communicated to me by R.~Bezrukavnikov.
\end{rem}

\Section{Conilpotent Comodules and Koszulity}

\subsection{Coaugmentation filtrations}
 Let $C$ be a coaugmented coalgebra with the coaugmentation map
$\c\:\k\rarrow C$.
 Then the {\it coaugmentation filtration\/}~$N$ on the coalgebra~$C$
is an increasing filtration defined by the formula~\cite{PV}
  $$
   N_nC=\{ c\in C \mid \D^{(n+1)}(c)\in C^{\ot n+1}_\c
   = \sum_{i=1}^{n+1} C^{\ot i-1}\ot\c(\k)\ot C^{\ot n-i+1}
   \sub C^{\ot n+1} \},
  $$
where $\D^{(m)}\:C\rarrow C^{\ot m}$ denotes the iterated
comultiplication map.
 In particular, we have $N_0C=\c(\k)$.

 Furthermore, let $P$ be a comodule over a coaugmented coalgebra $C$.
 Then the {\it coaugmentation filtration\/}~$N$ on the comodule~$P$
is defined as
 $$
  N_nP=\{ p\in P \mid \D'(p)\in N_nC\ot P\sub C\ot P \}.
 $$

 A coaugmented coalgebra $C$ is called {\it conilpotent\/}
if the coaugmentation filtration $N$ is full, that is $C=\bigcup_nN_nC$.
 Note that in this case for any $C$\+comodule $P$ one has
$P=\bigcup_nN_nP$; in other words,
{\it any comodule over a conilpotent coalgebra is conilpotent}.

 For any coaugmented coalgebra~$C$, we denote by $\Nilp C$
the maximal conilpotent subcoalgebra $\bigcup_nN_nC$
of the coalgebra~$C$.
 The following simple result is needed here.

\begin{thm}{Proposition 6}
 The coaugmentation filtrations respect the coalgebra and comodule
structures on a coalgebra~$C$ and comodule~$P$, that is
  $$
   \D(N_nC)\sub \sum_{i+j=n} N_iC\ot N_jC
   \qquad \text{and} \qquad
   \D'(N_nP)\sub \sum_{i+j=n} N_iC\ot N_jP.
  $$
 Furthermore, the associated graded coalgebra
$\gr_NC=\bop_{n=0}^\infty N_nC/N_{n-1}C$ and the comodule
$\gr_NP=\bop_{n=0}^\infty N_nP/N_{n-1}P$ over it
are one-cogenerated.
 \end{thm}

\pr{Proof}:
 The statements concerning the coalgebras only are proven in~\cite{PV}.
 Let us show that for any comultiplicative filtration~$N$
on a coalgebra~$C$ the filtration on a $C$\+comodule~$P$ given by
the above formula is compatible with the coaction.
 Let $\psi\:C\rarrow\k$ be a linear function annihilating $N_{k-1}C$;
we have to prove that $(\psi\ot\nbk\id)\D'(N_nP)\sub N_{n-k}P$.
 The latter inclusion, by the definition, can be rewritten as
$\D'(\psi\ot\nbk\id)\D'(N_nP)\sub N_{n-k}C\ot P$.
 Now we have
  \begin{align*}
    \D'(\psi\ot\nbk\id)\D'(N_nP)
      &=(\psi\ot\nbk\id\ot\nbk\id)(\id\ot\nbk\D')\D'(N_nP) \\*
      &=(\psi\ot\nbk\id\ot\nbk\id)(\D\ot\nbk\id)\D'(N_nP) \\*
      &\sub (\psi\ot\nbk\id\ot\nbk\id)(\D\ot\nbk\id)(N_nC\ot\nbk P)
      \,\sub\, N_{n-k}C\ot P.
  \end{align*}

 Since one has $\D'(p)\notin N_{n-1}C\ot P$ for any $p\notin N_{n-1}P$,
the last assertion of the proposition immediately follows.
 \qed

\subsection{Pro-finite groups}
 Let $G$ be a pro-finite group and $\k$ be a field.
 Consider the coalgebra~$\k(G)$ of locally constant $\k$-valued functions
on~$G$ with respect to the convolution; in other words,
let $\k(G)=\varinjlim\k(G/U)$, where the direct limit is taken over all
open normal subgroups $U\sub G$ and the finite-dimensional coalgebra
$\k(G/U)$ is the dual vector space to the group algebra $\k[G/U]$.
 Then the category of discrete $G$\+modules over~$\k$ is equivalent to
the category of $\k(G)$\+comodules.
 (Of course, the same is true over~$\Z$ or an arbitrary
commutative coefficient ring.)

 Now let $\c:\k\rarrow\k(G)$ be the coaugmentation map that takes
a constant from~$\k$ to the corresponding constant function on $G$.
 Then for any discrete $G$\+module $P$ over~$\k$ one has
$H^*(G,P)=H^*(\k(G),P)$, because $\Ext_G^*(\k,P)=\Ext_{\k(G)}^*(\k,P)$.

 Suppose $\k$~is a field of characteristic~$l$.
 Then one can see that the coaugmented coalgebra $\k(G)$ is conilpotent
if and only if $G$ is a pro-$l$-group.
 Moreover, for any pro-finite group~$G$ the coalgebra $\Nilp\k(G)$
is isomorphic to the group coalgebra $\k(G^{(l)})$
of the maximal quotient pro-$l$-group $G^{(l)}$ of~$G$.

\subsection{Nilpotency and Koszulity}
 The following theorem is the main result of my paper
with A.~Vishik~\cite{PV}.

\begin{thm}{Theorem 3}
 Let $A=H^*(C)$ be the cohomology algebra of a conilpotent
coalgebra~$C$.
 Assume that
  \begin{itemize}
    \item[(1)] the quadratic algebra\/ $\q A$ is Koszul;
    \item[(2)] the morphism of graded algebras\/ $\q A\rarrow A$
          is an isomorphism in degree~2 and a monomorphism
          in degree~3.
  \end{itemize}
 Then the algebra $A$ is quadratic (and therefore, Koszul).
 In addition, the graded coalgebra\/ $\gr_NC$ is Koszul and
there is an isomorphism $A\simeq (\gr_NC)^!$.
\qed
\end{thm}

 Our next theorem is a module version of the above one.
 Their proofs are completely parallel; still we prefer to give
the full argument here.

\begin{thm}{Theorem 4}
 Let $C$ be a conilpotent coalgebra such that the cohomology algebra
$A=H^*(C)$ is Koszul (equivalently, a coalgebra $C$ satisfies
the conditions of Theorem~3).
 Let $P$ be a comodule over~$C$.
 Consider the cohomology module $M=H^*(C,P)$ over the algebra~$A$.
 Assume that
  \begin{itemize}
    \item[(1)] the quadratic $A$\+module\/ $\q_AM$ is Koszul;
    \item[(2)] the morphism of graded $A$\+modules\/ $\q_AM\rarrow M$
          is an isomorphism in degree~1 and a monomorphism
          in degree~2.
  \end{itemize}
 Then the $A$\+module $M$ is quadratic (and therefore, Koszul).
 In addition, the graded\/ $\gr_NC$\+module\/ $\gr_NP$ is Koszul
and there is an isomorphism $M\simeq (\gr_NP)_{\gr_NC}^!$
of graded modules over the graded algebra $A\simeq(\gr_NC)^!$.
\end{thm}

\pr{Proof}:
 The coaugmentation filtrations on the coalgebra~$C$ and comodule~$P$
induce multiplicative filtrations on their cobar-complexes
by the standard rule
  $$
   N_nC^{+\ot i}=\sum_{j_1+\ds+j_i=n}N_{j_1}C^+\ot\ds\ot N_{j_i}C^+
  $$
and
  $$
   N_n(C^{+\ot i}\ot P)=\sum_{j_1+\ds+j_i+k=n}
   N_{j_1}C^+\ot\ds\ot N_{j_i}C^+\ot N_kP,
  $$
where $N_jC^+=N_jC/\c(\k)$; so the filtrations on the terms
$C^{+\ot i}$ and $C^{+\ot i}\ot P$ start with $N_i$.
 Clearly, the associated graded complexes coincide
with the cobar-complexes of $\gr_NC$ and $\gr_NP$.
 Therefore, we obtain a multiplicative spectral sequence
  $$
   E_1^{i,j}=H^{i,j}(\gr_NC) \implies H^i(C)
  $$
and a module spectral sequence
  $$
   \Ee_1^{i,j}=H^{i,j}(\gr_NC,\gr_NP) \implies H^i(C,P)
  $$
over it.
 Since the filtrations on the cobar-complexes are increasing ones,
these spectral sequences converge in the direct limit.
 The differentials have the form
$d_r\:E_r^{i,j}\rarrow E_r^{i+1,j-r}$
and $\dd_r\:\Ee_r^{i,j}\rarrow\Ee_r^{i+1,j-r}$.
 Furthermore, there are induced increasing filtrations $N$ on
$H^*(C)$ and $H^*(C,P)$ compatible with the multiplication
and the action and such that $\gr_N^jH^i(C)=E_\infty^{i,j}$
and $\gr_N^jH^i(C,P)=\Ee_\infty^{i,j}$.

 It was shown in the proof of Theorem~3 given in~\cite{PV}
that the spectral sequence $E_r^{i,j}$ degenerates at~$E_1$
with $E_1^{i,j}=E_\infty^{i,j}=0$ for $i\ne j$
and $E_1^{i,i}=E_\infty^{i,i}=H^i(C)=A_i$.
 In particular, we have $N_iH^i(C)=H^i(C)$.
 Besides, by Theorem~3 the graded algebra $\gr_NC$ is Koszul
and $H^*(\gr_NC)\simeq(\gr_NC)^!\simeq A$.

 {}From the spectral sequence $\Ee_r^{i,j}$ we see that the submodule
$\bop_{i=0}^\infty N_iH^i(C,P)$ of the $A$\+module $M=H^*(C,P)$
is isomorphic to the quotient module of the diagonal cohomology
$A$\+module $\bop_{i=0}^\infty H^{i,i}(\gr_NC,\gr_NP)$
by the images of the differentials.
 According to Proposition~6, the graded $\gr_NC$\+comodule $\gr_NP$
is one-cogenerated, hence (by Proposition~2)
we have $\Ee_1^{0,j}=H^{0,j}(\gr_NP,\gr_NC)=0$ for $j>0$,
which implies $H^0(C,P)=N_0H^0(C,P)\simeq H^{0,0}(\gr_NC,\gr_NP)$.
 Now the spectral sequence shows that $E_\infty^{1,1}=E_1^{1,1}$
and therefore $N_1H^1(C,P)\simeq H^{1,1}(\gr_NC,\gr_NP)$.
 Since (by Proposition~3) the diagonal cohomology $A$\+module
$\bop_{i=0}^\infty H^{i,i}(\gr_NC,\gr_NP)$ is quadratic,
we conclude that it is isomorphic to~$\q_A H^*(C,P)$.

 By Proposition~3 again, the $A$\+module
$\bop_{i=0}^\infty H^{i,i}(\gr_NC,\gr_NP)$
is quadratic dual to the $\gr_NC$\+comodule $\q_{\gr_NC}\gr_NP$;
since we suppose that the $A$\+module $\q_A H^*(C,P)$ is Koszul,
the dual ${\gr_NC}$\+comodule $\q_{\gr_NC}\gr_NP$ is Koszul, too
(Proposition~4).
 Besides, we have assumed that the morphism
$\q_AH^*(C,P)\rarrow H^*(C,P)$ is an isomorphism in degree~$1$,
hence $H^1(C,P)=N_1H^1(C,P)$ and $\Ee_\infty^{1,n}=0$.
 Furthermore, we have supposed that the morphism
$\q_AH^*(C,P)\rarrow H^*(C,P)$ is a monomorphism in degree~$2$.
 Since $\Ee_1^{2,2}=H^{2,2}(\gr_NC,\gr_NP)\simeq\q_AH^*(C,P)_2$
and $\Ee_\infty^{2,2}\simeq N_2H^2(C,P)$,
it follows that the map $\Ee_1^{2,2}\rarrow\Ee_\infty^{2,2}$
is a monomorphism and therefore all the differentials
$\dd_r\:\Ee_r^{1,2+r}\rarrow \Ee_r^{2,2}$
targeting at $\Ee_r^{2,2}$ vanish.

 Now let us establish by induction on~$j$ that
$H^{1,j}(\gr_NC,\gr_NP)=0$ for all $j>1$.
 Assume that this is true for $1<j\le n-1$.
 By Proposition~2, it follows that the map
$r_{gr_NC,\,\gr_NP}\:\gr_NP\rarrow\q_{\gr_NC}\gr_NP$
is an isomorphism in degree $\.\le n-1$.
 Therefore, the induced map of the cobar-complexes is also
an isomorphism in these degrees, hence in particular
$H^{1,j}(\gr_NC,\gr_NP)=H^{1,j}(\gr_NC,\q_{\gr_NC}\gr_NP)$
for $j\le n-1$ (and actually even for $j\le n$).
 Since the $\gr_NC$\+comodule $\q_{\gr_NC}\gr_NP$ is Koszul,
it follows that $\Ee_1^{2,j}=H^{2,j}(\gr_NC,\gr_NP)=0$
for all $2<j\le n-1$.
 From the latter it is clear that $\Ee_1^{1,n}=\Ee_\infty^{1,n}$,
since the differentials targeting at $\Ee_r^{2,2}$ also vanish.
 It was noticed above that $\Ee_\infty^{1,n}=0$;
so $H^{1,n}(\gr_NC,\gr_NP)=\Ee_1^{1,n}=0$ and we are done.

 We have shown that the $\gr_NC$\+comodule $\gr_NP$ is quadratic
and the comodule $\q_{\gr_NC}\gr_NP$ is Koszul,
hence the comodule $\gr_NP$ is Koszul.
 It follows that $\Ee_1^{i,j}=0$ for $i\ne j$, so the spectral sequence
degenerates and $H^*(C,P)=H^*(\gr_NC,\gr_NP)=(\gr_NP)_{\gr_NC}^!$.
 Therefore, $H^*(C,P)$ is a Koszul $H^*(C)$\+module, too.
\qed

\Section{Maximal Conilpotent Subcoalgebra and Koszulity}

 The results of this section connect the cohomology of a coaugmented
coalgebra with that of its maximal conilpotent subcoalgebra
(and, in particular, the cohomology of a pro-finite group with
the cohomology of its maximal quotient pro-$l$-group)
under certain Koszulity conditions.
 Let us start with the following important lemma.

\begin{thm}{Lemma 1}
 Let $C$ be a coaugmented coalgebra and\/ $\Nilp C$ be its maximal 
conilpotent subcoalgebra; then the natural map
$H^*(\Nilp C)\rarrow H^*(C)$ is an isomorphism in degree~1
and a monomorphism in degree~2.
 In particular, the quadratic algebras $\q H^*(C)$ and
$\q H^*(\Nilp C)$ are naturally isomorphic to each other.
\end{thm}

\pr{Proof}:
 Any comodule over the coalgebra $\Nilp C$ can be considered
as a comodule over~$C$ in the most obvious way.
 Furthermore, a $C$\+comodule comes from a comodule over $\Nilp C$
if and only if it is a direct limit of successive extensions
of the trivial $C$\+comodule~$\k$ (which is defined in terms
of the coaugmentation).
 This is not difficult to prove, in the ``if'' direction by induction,
using the coaugmentation filtrations on comodules and Proposition~6.
 It is helpful to remember that every comodule over a coassociative
coalgebra is a direct limit of finite-dimensional comodules.

 For any two comodules $P$ and~$Q$ over the subcoalgebra $\Nilp C$
there is a natural morphism of graded vector spaces
$\Ext^*_{\Nilp C}(P,Q)\rarrow\Ext^*_C(P,Q)$.
 Now it follows from the above that this map is always
an isomorphism in degree~1.
 It turns out that one can derive purely formally from the latter
that this is also a monomorphism in degree~2.
 Namely, it suffices to compute the spaces~$\Ext^2$ in question
in terms of~$\Ext^1$ using a one-step injective resolution
of the comodule~$Q$ over the coalgebra $\Nilp C$, for example
$Q\rarrow \Nilp C\ot Q\rarrow (\Nilp C\ot Q)/Q$.
\qed

\smallskip

 The next result can be considered as a direct generalization
of Theorem~3.

\begin{thm}{Theorem 5}
 Let $C$ be a coaugmented coalgebra such that the cohomology
algebra $H^*(C)$ satisfies the conditions of Theorem~3.
 Then the cohomology algebra of the maximal conilpotent subcoalgebra 
$\Nilp C$ is isomorphic to the quadratic part of the cohomology
of~$C$, that is $H^*(\Nilp C)\simeq \q H^*(C)$.
\end{thm}

\pr{Proof}:
 Let us show that the cohomology algebra of the coalgebra $\Nilp C$
satisfies the conditions of Theorem~3, too.
 Indeed, the second statement of Lemma~1 provides the condition~(1).
 To check the condition~(2), one should consider the composition
$\q H^*(C)\simeq\q H^*(\Nilp C)\rarrow H^*(\Nilp C)\rarrow H^*(C)$
and use the first statement of Lemma~1.
 Applying Theorem~3, we conclude that the algebra $H^*(\Nilp C)$
is quadratic. 
 The latter implies the desired isomorphism of graded algebras.
\qed

\begin{thm}{Corollary 1}
 Let $C$ be a coaugmented coalgebra and $\Nilp C$ be its maximal
conilpotent subcoalgebra.
 In this setting, if the cohomology algebra $H^*(C)$ is Koszul,
then the natural homomorphism $H^*(\Nilp C)\rarrow H^*(C)$
is an isomorphism.  \qed
\end{thm}

 The analogous results comparing the cohomology of a pro-finite group
with that of its maximal quotient pro-$l$-group can be either deduced
from the results for coalgebras using the observations made in
subsection~3.2, or they can be proven independently in the very
parallel way.
 The proof of the analogue of Lemma~1 can be done differently
(and even simpler) for pro-finite groups.

\begin{thm}{Lemma 2}
 If $G$ is a pro-finite group and $G^{(l)}$ is its maximal quotient
pro-$l$-group, then the natural map
$H^i(G^{(l)},\.\Z/l)\rarrow H^i(G,\.\Z/l)$ is an isomorphism for $i=1$
and a monomorphism for $i=2$.
 In particular, the quadratic algebras $\q H^*(G,\.\Z/l)$ and
$\q H^*(G^{(l)},\.\Z/l)$ are naturally isomorphic to each other.
\end{thm}

\pr{Proof}:
 Let $G^{(l)}=G/G'$; then by the definition the subgroup $G'$ has
no nontrivial quotient pro-$l$-groups, so $H^1(G',\.\Z/l)=0$.
 It remains to apply the Serre--Hochschild spectral sequence
$E_2^{p,q}=H^p(G^{(l)},H^q(G',\.\Z/l))\implies H^{p+q}(G,\.\Z/l)$.
\qed

\begin{thm}{Corollary 2}
 Let $G$ be a pro-finite group and $G^{(l)}$ be its maximal quotient
pro-$l$-group.
 If the cohomology algebra $H^*(G,\.\Z/l)$ satisfies the conditions
of Theorem~3, then the cohomology algebra of the group~$G^{(l)}$
is isomorphic to the quadratic part of the cohomology of~$G$,
 $$
   H^*(G^{(l)},\.\Z/l)\,\simeq\, \q H^*(G,\.\Z/l).
 $$
 Furthermore, if the cohomology algebra $H^*(G,\.\Z/l)$ is Koszul,
then the natural homomorphism
 $$
   H^*(G^{(l)},\.\Z/l)\lrarrow H^*(G,\.\Z/l)
 $$
is an isomorphism. \qed
\end{thm}

\begin{rem}{Remark:}
 All the results of this section actually have a rather general
categorical nature: they can be formulated for an arbitrary
abelian category (or, better yet, a triangulated category)
in place of the category of $C$\+comodules with a ``nilpotent''
abelian subcategory in place of the subcategory of comodules
over $\Nilp C$; the nilpotency can be generalized to mean
that all the objects have finite length.
 The proof is based on the generalization of Theorem~3
to nilpotent abelian categories.
 We wouldn't go so far here, but will only mention several
elementary applications.

 Let $A$ be an augmented algebra.
 The {\it coalgebra of pro-nilpotent completion\/} $C=A^\wedge$
of the algebra~$A$ is defined as $C=\varinjlim_I (A/I)^*$,
where the limit goes over all the ideals $I\sub A_+$ such that
the algebra $A/I$ is a finite-dimensional nilpotent algebra
(meaning that there exists such~$n$ that $(A_+/I)^n=0$)
and $(A/I)^*$ is the dual vector space to $A/I$ with its natural
coalgebra structure.
 Then the following analogue of Theorem~5 is true: if the
cohomology algebra $H^*(A)=\Ext^*_A(\k,\k)$ satisfies the 
conditions of Theorem~3, then the algebra $H^*(C)$ is
isomorphic to the quadratic part of~$H^*(A)$.

 Furthermore, the analogue of Corollary~2 is valid for a discrete
group $\G$ in place of~$G$ and its pro-$l$-completion $\G^{(l)}$
or pro-unipotent completion $\G^{(\Q)}$ in place of~$G^{(l)}$,
where one considers the cohomology with the constant coefficients
$\Z/l$ or~$\Q$, respectively.
 One can reduce the last two statements to the one about
the cohomology of an augmented algebra using the results
of Quillen's paper~\cite{Quil}, or prove all the three
assertions independently just in the same way as it was
done for coalgebras above.
\end{rem}

\Section{Morphisms of Graded Algebras and Koszulity}

 The results of this and the next sections were found in an attempt
to generalize and clarify some of the statements from the paper
of J.~Backelin and R.~Fr\"oberg~\cite{BF}.

 Suppose we are given a homomorphism of algebras $f\:A\rarrow B$.
 Then there is an induced structure of left $A$\+module on
the vector space~$B$.
 For any right $A$\+module $M$ and left $B$\+module $N$
there is a natural isomorphism $(M\ot_AB)\ot_BN\simeq M\ot_AN$,
which implies the analogous isomorphism on the level of
derived functors.
 Thus we obtain a spectral sequence
$$
 E^2_{p,q}=\Tor_p^B(\Tor_q^A(M,B),N)\implies \Tor_{p+q}^A(M,N)
$$
with the differentials $d^r\:E^r_{p,q}\rarrow E^r_{p-r,q+r-1}$.
 In particular, for a homomorphism of augmented algebras~$f$
and the trivial modules $M=\k$ and $N=\k$ we get
$$
 E^2_{p,q}=H_p(B^\opp,H_q(A,B))\implies H_{p+q}(A).
$$
 Recall that in our notation $\Tor^A_*(M,\k)=H_*(A^\opp,M)$ and
$\Tor^B_*(\k,N)=H_*(B,N)$; besides,
$H_*(A)=H_*(A^\opp)=\Tor_*^A(\k,\k)$.

 If $A$ and~$B$ are graded algebras ($M$ and $N$ are graded modules)
and the homomorphism~$f$ preserves the grading, then all the terms
of the spectral sequence $E^r_{p,q}$ bear the corresponding additional
grading which is preserved by all the differentials.

 We will be interested in two specific cases, which can be
thought of as {\it two kinds of Koszulity condition
on a morphism\/ $f\:A\rarrow B$}.
 The first case is when one has $H_{i,j}(A,B)=0$ unless $i=j$.
 A bit more general, but nevertheless interesting, second case is when
$H_{i,j}(A,B)=0$ unless $j-i\le 1$.
 Note that the first condition implies surjectivity of the morphism~$f$;
the second one doesn't.

 For a nonnegatively graded vector space $N=\bop_{i=0}^\infty N_i$
and an integer~$n\ge0$ let us denote by $N(n)$ the graded vector
space with the components $N(n)_i=N_{i-n}$.
 We will say that a graded module $M$ over a Koszul algebra $A$
is a {\it Koszul module in the grading shifted by\/~$n$} if there
exists a Koszul $A$\+module~$N$ such that $M=N(n)$.

\begin{thm}{Theorem 6}
 Let $f\:A\rarrow B$ be a homomorphism of graded algebras
such that the corresponding vector spaces $H_{i,j}(A,B)$
are zero unless $j-i\le 1$.
 Then if the graded algebra $A$ is Koszul, then the algebra
$B$ is Koszul also.
 Furthermore, if $H_{i,j}(A,B)=0$ unless $i=j$, then
the graded algebras $A$ and $B$ are Koszul simultaneously.
\end{thm}

\pr{Proof}:
 We have to prove that $H_{p,j}(B)=0$ for all $j\ne p$.
 Proceeding by induction, assume that this is true for all
$0\le p<n$.
 Since $H_{q,j}(A,B)=0$ for all $q\ne j$,~$j+1$ and 
the graded right $B$\+module $H_q(A,B)$ can be presented as
an extension of its grading components endowed
with the trivial $B$\+module structures,
it follows immediately that for any $p<n$ the term
$E^2_{p,q}=H_p(B^\opp,H_q(A,B))$ of the above spectral sequence
has no components of degrees other than $p+q$ and $p+q+1$.

 Now let us consider the term $E^2_{n,0}$.
 The graded components of the terms $E^r_{n,0}$ can only cancel
in the spectral sequence with those of the terms $E^r_{n-r,r-1}$;
according to the above, the latter are concentrated in the degrees
$n-1$ and~$n$.
 The limit term $E^\infty_{n,0}=\gr_n H_n(A)$ is known to be
only nonzero in degree~$n$; it follows easily that the graded vector
space $H_n(B,H_0(A,B))=E^2_{n,0}$ is concentrated in degree~$n$.

 By assumption, the graded vector space $H_0(A,B)=B/f(A)$
is concentrated in the degrees $0$ and~$1$.
 Thus we have an exact triple of graded right $B$\+modules
  $$
   0\lrarrow B_1/f(A_1)\.(1)\lrarrow H_0(A,B)\lrarrow \k\lrarrow0,
  $$
where the vector spaces $B_1/f(A_1)$ and~$\k$ have the trivial
$B$\+module structures.
 Consider the corresponding long exact sequence of homology:
the term $H_n(B)$ is placed between $H_n(B^\opp,H_0(A,B))$
and $H_{n-1}(B)\ot B_1/f(A_1)\.(1)$ and so it follows from
the induction hypothesis that $H_{n,j}(B)=0$ for all $j\ne n$.

 Conversely, if the graded spaces $H_p(B)$ are concentrated
in degree~$p$ and the spaces $H_q(A,B)$ are concentrated
in degree~$q$, then all the terms $E^2_{p,q}=H_p(B^\opp,H_q(A,B))$
are concentrated in degree $p+q$.
 In this case the same holds for the terms
$E^\infty_{p,q}=\gr_p H_{p+q}(A)$, hence 
the spaces $H_n(A)$ are concentrated in degree~$n$.
\qed

\smallskip

 The next corollary lists the most typical specific situations
(see also~\cite{BF}).

\begin{thm}{Corollary 3}
 Let $f\:A\rarrow B$ be a homomorphism of graded algebras.
 Suppose that the algebra $A$ is Koszul and one of the following
three conditions holds:
  \begin{itemize}
    \item[(a)] the $A$-module $B$ is Koszul (equivalently,
               the morphism $f$ is surjective and its kernel $J$
               is a Koszul $A$-module in the grading shifted by~$1$), or
    \item[(b)] the morphism $f$ is injective and its cokernel $B/f(A)$
               is a Koszul $A$-module in the grading shifted by~$1$, or
    \item[(c)] the morphism $f$ is surjective and its kernel $J$
               is a Koszul $A$-module in the grading shifted by~$2$.
  \end{itemize}
 Then the algebra $B$ is Koszul, too.  \qed
\end{thm}

 The proof of the following result is completely analogous
to that of Theorem~6.

\begin{thm}{Proposition 7}
 Let $f\:A\rarrow B$ be a morphism of graded algebras.
 Assume that $H_{0,j}(A,B)=0$ for all $j\ne 0$,~$1$ and
$H_{i',\,j}(B^\opp,H_{i''}(A,B))=0$ for all $i''\ge1$ and
$j-i'-i''\ne 0$,~$1$.
 In this case if the algebra~$A$ is Koszul, then the algebra~$B$
is Koszul, too.
 If the algebra~$B$ is Koszul and all the right $B$\+modules
$H_i(A,B)$ are Koszul in the grading shifted by~$i$, then
the algebra~$A$ is Koszul.  \qed
\end{thm}

\begin{thm}{Corollary 4}
 Let $f\:A\rarrow B$ be a homomorphism of graded algebras.
 Assume that $A$ is a Koszul algebra and one of the following
three conditions holds:
  \begin{itemize}
    \item[(a)] the algebra~$B$ is a free left $A$\+module and
               the right $B$\+module $H_0(A,B)=B/A_+B$ is Koszul, or
    \item[(b)] the morphism~$f$ is surjective, its kernel~$J$
               is a free left $A$-module, and the right $B$\+module
               $H_0(A,J)=H_1(A,B)$ is Koszul in the grading shifted
               by~$1$, or
    \item[(c)] same as (b), except that the module
               is Koszul in the grading shifted by~$2$.
  \end{itemize}
 Then the algebra~$B$ is Koszul, too.  \qed
\end{thm}

\begin{thm}{Proposition 8}
 Let $f\:A\rarrow B$ be a morphism of Koszul algebras.
 In this case if the algebra~$B$ is a free left $A$\+module,
then the right $B$\+module $H_0(A,B)$ is Koszul.
 If the morphism~$f$ is surjective and its kernel~$J$ is a free
left $A$\+module, then one has $H_{i,j}(B^\opp, H_0(A,J))=0$
for all $j-i\ne 1$,~$2$.
 More specifically, the right $B$\+module $H_0(A,J)$ is Koszul
in the grading shifted by~$1$ iff the morphism of dual coalgebras
$f^?\:H_*(A)\rarrow H_*(B)$ is surjective;
the module $H_0(A,J)$ is Koszul in the grading shifted
by~$2$ iff the morphism~$f^?$ is injective.
\end{thm}

\pr{Proof:}
 All the statements follow from the spectral sequence
of the morphism~$f$.
\qed

\Section{Morphisms of Koszul Algebras: Koszulity and Freeness}

 Let $f\:A\rarrow B$ be a homomorphism of Koszul algebras.
 Note that just as the morphism~$f$ makes the algebra~$B$ a left
and a right $A$\+module, the dual morphism $f^?\:A^?\rarrow B^?$
makes the coalgebra~$A^?$ a left and a right $B^?$\+comodule.
 In this section we describe the correspondence between the
homological conditions of the above type on the morphisms
$f$ and~$f^?$.
 Here is the principal result.

\begin{thm}{Theorem 7}
 Let $f\:A\rarrow B$ be a morphism of Koszul algebras and
$f^?\:A^?\rarrow B^?$ be the dual morphism of Koszul coalgebras.
 Then there are natural isomorphisms of the (co)module (co)homology
vector spaces
 $$
  H_{i,j}(A,B)\,\simeq\,H^{j-i,\,j}(B^{?\opp},\.A^{?\opp})
 $$
compatible with the right action of~$B$ and the left coaction
of~$A^?$ on both sides.
\end{thm}

\pr{Proof}:
 According to Proposition~4, the Koszul complex $K(A^?,A)$
with the homological grading $K_i(A^?,A)=A^?_i\ot A$ is
a free graded right $A$\+module resolution of the trivial
$A$\+module~$\k$.
 Therefore, the space $H_{i,j}(A,B)$ can be computed as
the homology space of the complex $K(A^?,A)\ot_AB=A^?\ot B$
at the term $A^?_i\ot B_{j-i}$.
 Analogously, the Koszul complex $K(B^?,B)$ with the
cohomological grading $K^i(B^?,B)=B^?\ot B_i$ is a cofree
graded left $B^?$\+comodule resolution of the trivial
$B^?$\+comodule~$\k$.
 It follows that the space $H^{i,j}(B^{?\opp},\.A^{?\opp})$
is isomorphic to the cohomology space of the complex
$A^?\oa_{B^?}K(B^?,B)=A^?\ot B$ at the term $A^?_{j-i}\ot B_i$.
 Now we see that both sides of the desired isomorphism
of vector spaces are computed by the Koszul complex
$K(f;\.A^?,B)=A^?\ot B$ with the differential given by
the formula
 \begin{gather*}
   \d=(\id\ot m)(\id\ot f_1\ot\id)(\D\ot\id)\:
    A^?\ot B\lrarrow A^?\ot A_1\ot B  \\
    \lrarrow A^?\ot B_1\ot B\lrarrow A^?\ot B.
 \end{gather*}
 The right action of~$B$ and the left coaction of~$A^?$ on
both $H_*(A,B)$ and $H^*(B^{?\opp},A^{?\opp})$ come from
the natural action and coaction on this Koszul complex.
\qed

\begin{thm}{Corollary 5}
 The conditions (a),~(b),~(c) of Corollaries~3 and~4 are
respectively dual to each other in the following sense.
 Let $f\:A\rarrow B$ be a morphism of Koszul algebras and
$f^?\:A^?\rarrow B^?$ be the dual morphism of Koszul coalgebras.
 Then the morphism $f$ satisfies the right module version of
a condition from Corollary~3 if and only if the morphism~$f^?$
satisfies the coalgebra version of the corresponding condition from
Corollary~4, and vice versa.
 Specifically, the following statements hold.
  \begin{itemize}
    \item[(a)] The algebra $B$ is a free right $A$\+module if
               and only if the coalgebra~$A^?$ is a Koszul
               left $B^?$\+comodule.
               If this is the case, the left $B$\+module $H_0(A^\opp,B)$
               is the Koszul module dual to the left $B$\+comodule~$A^?$.
  \item[(b-c)] If the morphism~$f$ is surjective, then its kernel~$J$
               is a free right $A$\+module if and only if
               $H^{i,j}(B^?,A^?)=0$ for all $j-i\ne 0$,~$1$.
               (One has $H_{i,i}(B^?,A^?)=0$ for all $i>0$
               and $H_{j-1,\,j}(B^?,A^?)\simeq H_{0,j}(A^\opp,J)$
               for all~$j$.)
  \end{itemize}
 In the second case, if~(b) the morphism~$f^?$ is surjective,
then its kernel~$K$ and the space $H_0(A^\opp,J)$ of graded
generators of the right ideal~$J$ are dual left Koszul
$B^?$\+comodule and $B$\+module, both in the grading shifted by~1.
 If~(c) the map~$f_1$ is an isomorphism, then the cokernel~$C$
of the morphism~$f$ and the generator space $H_0(A^\opp,J)$ are
dual left Koszul $B^?$\+comodule and $B$\+module, both
in the grading shifted by~2. \qed
\end{thm}

\Section{Graded Lie Coalgebras and Koszulity}

 It turns out somewhat unexpectedly that the results Theorem~7 implies
for morphisms of Koszul commutative and Lie (co)algebras are stronger
than in general.
 On the other hand, Lie coalgebras involve some troubles~\cite{Mich},
which disappear in the positively graded case.

 By a {\it graded Lie algebra\/} we mean here a positively graded
vector space $L=\bop_{i=1}^\infty L_n$ equipped with a Lie algebra
structure such that $[L_i,L_j]\sub L_{i+j}$ for all $i$,~$j$.
 The universal enveloping algebra $U(L)$ of a graded Lie algebra
is a graded associative algebra.
 A graded Lie algebra is called {\it quadratic\/} if it isomorphic
to the quotient algebra of the free Lie algebra generated by~$L_1$
by an ideal generated by a set of elements of degree~$2$.
 A graded Lie algebra~$L$ is quadratic if and only if the algebra
$U(L)$ is quadratic.
 The homology coalgebra of a Lie algebra~$L$ can be defined as
the homology of its enveloping algebra, $H_*(L)=H_*(U(L))$, or,
more explicitly, as the homology coalgebra of the standard complex
 $$
  H_*(L)=H_*(\.\k\larrow L\larrow \bwe^2 L\larrow \bwe^3 L\larrow\dsb).
 $$
 From the second definition it is obvious that the homology coalgebra
of a Lie algebra is always skew-cocommutative.
 By a {\it graded super-Lie algebra\/} we mean a positively graded
vector space equipped with a super-Lie algebra structure with
respect to the parity induced by the grading such that the bracket
respects the grading.
 There is a similar standard complex computing the homology of
a super-Lie algebra.

\begin{thm}{Proposition 9}
 A quadratic algebra~$A$ is isomorphic to the enveloping algebra
of a quadratic Lie algebra if and only if the dual quadratic
coalgebra~$A^?$ is skew-cocommutative.
 A quadratic algebra~$A$ is isomorphic to the enveloping algebra
of a quadratic super-Lie algebra iff the coalgebra~$A^?$
is cocommutative.  \qed
\end{thm}

 A graded Lie algebra or super-Lie algebra~$L$ is called {\it Koszul\/}
if one has $H_{i,j}(L)=0$ for $i\ne j$; in other words, $L$~is Koszul
if the enveloping associative algebra~$U(L)$ is Koszul.
 Let us denote by $L^?=U(L)^?$ the quadratic coalgebra dual to $U(L)$.

\begin{thm}{Corollary 6}
 Let $\phi\:L'\rarrow L''$ be a morphism of Koszul Lie algebras
or super-Lie algebras.
 Then the coalgebra $L^{\1\,?}$ is a Koszul $L^{\1\1\,?}$\+comodule
if and only if the morphism~$\phi$ is injective.
 Furthermore, one has $H^{i,j}(L^{\1\1\,?},\.L^{\1\,?})=0$
for all $j-i\ne 0$,~$1$ if and only if the kernel~$L$
of the morphism~$\phi$ is a free graded (super)\+Lie algebra.
\end{thm}

\pr{Proof}:
 It follows from the Poincar\'e--Birkhoff--Witt theorem that
the enveloping algebra $U(L'')$ is a free graded $U(L)$\+module
if and only if the morphism~$\phi$ is injective, so it remains
to apply Corollary~5\,(a) to prove the first statement.
 For an arbitrary morphism of Lie algebras~$\phi$, it is clear
that the algebra $U(L'')$ is a free module over the algebra
$U(\im \phi)=U(L'/L)$, the $U(L')$\+module $U(L'/L)$ is induced
from the trivial $U(L)$\+module, and thus there are natural
isomorphisms $H_*(U(L'),\.U(L''))\simeq
H_*(U(L'),\.U(L'/L))\ot_{U(L'/L)}U(L'')
\simeq H_*(L)\ot_{U(L'/L)}U(L'')$.
 Therefore, one has $H_i(U(L'),\.U(L''))=0$ for $i\ge 2$ if and
only if~$L$ is a free Lie algebra and the second statement
follows from Theorem~7.
\qed

\smallskip

 Recall that a {\it Lie coalgebra\/} is a vector space $\La$ together
with a cobracket map $\de\:\La\rarrow\bwe^2\La$ satisfying the equation
dual to the Jacobi identity.
 The cohomology algebra of a Lie coalgebra $\La$ is defined as
the cohomology of the standard complex
 $$
  H^*(\La)=H^*(\.\k\rarrow \La\rarrow \bwe^2\La\rarrow
  \bwe^3\La\rarrow\dsb).
 $$
 The notion of the coenveloping coalgebra of a Lie coalgebra involves
troubles, which may be overcome by restricting to the following
condition.
 A Lie coalgebra is called {\it conilpotent} if it is a direct limit
of finite-dimensional Lie coalgebras dual to nilpotent Lie algebras.
 The {\it conilpotent enveloping coalgebra\/} $C(\La)$ of a conilpotent
Lie coalgebra~$\La$ is defined as the universal object in
the category of morphisms of conilpotent Lie coalgebras $C\rarrow\La$
from the Lie coalgebra associated with a conilpotent coassociative
coalgebra~$C$ to the Lie coalgebra~$\La$.
 One can show that for any conilpotent Lie coalgebra~$\La$ the cohomology
algebras $H^*(\La)$ and $H^*(C(\La))$ are naturally isomorphic.

 By a {\it graded Lie coalgebra\/} we mean a positively graded vector
space $\La=\bop_{n=1}^\infty\La_n$ endowed with a Lie coalgebra
structure such that $\de(\La_n)\sub\sum_{i+j=n}\La_i\wedge\La_j$.
 Obviously, any positively graded Lie coalgebra is conilpotent.
 The conilpotent coenveloping coalgebra of a graded Lie coalgebra
is naturally a graded coassociative coalgebra and coincides
with the coenveloping coalgebra defined with respect to the category
of graded coalgebras.
 There are analogous notions of a {\it conilpotent super-Lie coalgebra},
its coenveloping coalgebra, and its cohomology algebra.
 A graded (super)\+Lie coalgebra $\La$ is called {\it Koszul\/}
if one has $H^{i,j}(\La)=0$ for all $i\ne j$.

\begin{thm}{Corollary 7}
 Let $f\:A\rarrow B$ be a morphism of commutative
or skew-commutative Koszul algebras.
 Then the algebra~$B$ is a Koszul $A$\+module if and only if
the dual morphism $f^?\:A^?\rarrow B^?$ is surjective.
\end{thm}

\pr{Proof}:
 This follows from the appropriate version of Proposition~9
and the analogue of Corollary~6 for morphisms of Koszul
Lie or super-Lie coalgebras.
 An important difference is that the computation with enveloping
algebras from the proof of the latter Corollary wouldn't work for
coenveloping coalgebras  in general, since
the dual version of Poincare--Birkhoff--Witt theorem does not
always hold.
 It holds for conilpotent Lie coalgebras, though.
\qed

\smallskip

 By the definition, a {\it pro-Lie algebra\/} is a projective limit of
finite-dimensional Lie algebras and a {\it pro-nilpotent Lie algebra\/}
is a projective limit of finite-dimensional nilpotent Lie algebras.
 Note that the category of pro-Lie algebras is anti-equivalent
to the category of locally finite Lie coalgebras and the category of
pro-nilpotent Lie algebras is anti-equivalent to the category of
conilpotent Lie coalgebras.
 Any positively graded pro-Lie algebras is, of course, pro-nilpotent.
 The cohomology algebra of a pro-nilpotent Lie algebra is defined as
the cohomology of the corresponding Lie coalgebra; a graded
pro-nilpotent (super)\+Lie algebra is called Koszul if the corresponding
(super)\+Lie coalgebra is Koszul.
 A graded pro-nilpotent (super)\+Lie algebra is free iff
the corresponding (super)\+Lie coalgebra is conilpotent cofree.

\begin{thm}{Corollary 8}
 Let $\phi\:\L'\rarrow\L''$ be a morphism of Koszul pro-Lie algebras or
pro-super-Lie algebras.
 Then one has $H_{i,j}(H^*(\L''),\.H^*(\L'))=0$ for all $j-i\ne 0,1$
if and only if the kernel $\L$ of the morphism $\phi$ is a free
pro-(super)\+Lie algebra. 
\end{thm}

This is the dual version of the second statement of Corollary~6.
\qed

 The result about freeness of a graded pro-Lie algebra mentioned
in subsection~1.8 of the Introduction now follows immediately from
Theorem~6 and Corollary~8.

\Section{Main Theorem}

 Two different proofs of Theorem~2 are given below, one based
on Theorem~4 and the other on Theorem~7.
 Both approaches actually lead to somewhat stronger results, each
in its own direction; they are formulated here as Theorem~9 and
Corollary~9.

\begin{thm}{Theorem 8}
 Let $g\:C\rarrow D$ be a morphism of conilpotent coalgebras.
 Assume that the cohomology algebras $A=H^*(C)$ and $B=H^*(D)$
are Koszul.
 In this case, if for a certain integer $t\ge 1$ one has
$H_{i,j}(A^\opp,\.B^\opp)=0$ for all $j-i\ge t$ (or even only
for all $j-i=t$), then $H^i(D,C)=0$ for all $i\ge t$.
 For example, if the right $A$\+module $B$ is Koszul, then
the left $D$\+comodule $C$ is cofree.
 On the other hand, if the algebra~$B$ is Koszul, then the morphism
$C\rarrow D$ is injective and at the same time the $B$\+module
$H^*(D,C)$ is Koszul if and only if the algebra~$A$ is Koszul and
at the same time the right $A$\+module $B$ is free.
\end{thm}

\pr{Proof}:
 Note that any morphism of conilpotent coalgebras $g\:C\rarrow D$
commutes with the coaugmentations.
 For the coaugmentation filtrations~$N$ on the coalgebras $C$ and~$D$
one obviously has $g(N_iC)\sub N_iD$ and it follows that these
filtrations are compatible with the $D$\+comodule structure on~$C$.
 Thus there is a spectral sequence
 $$
  E_1^{i,j}=H^{i,j}(\gr_ND,\.\gr_NC)\implies H^i(D,C),
 $$
where the structure of graded $\gr_ND$\+comodule on $\gr_NC$
coincides with the one induced by the morphism of graded coalgebras
$\gr_Ng\:\gr_NC\rarrow \gr_ND$.
 Therefore, we have $H^t(D,C)=0$ provided that $H^t(\gr_ND,\.\gr_NC)=0$.
 According to Theorem~3, there are natural isomorphisms
$\gr_NC\simeq A^?$ and $\gr_ND\simeq B^?$.
 Now the first statement follows from Theorem~7; it only remains
to notice that for any comodule~$P$ over a conilpotent coalgebra~$C$
one has $H^i(C,P)=0$ for all $i\ge t\/$ whenever $H^t(C,P)=0$.
 The latter property is deduced from the fact that $H^0(C,P)=0$
implies $P=0$.

 In general, the coaugmentation filtration on the coalgebra~$C$ can
differ from the coaugmentation filtration defined on the space~$C$
as a $D$\+comodule.
 One can see that the two filtrations coincide if and only if
the morphism $C\rarrow D$ is injective, or, equivalently,
the morphism $A_1\rarrow B_1$ is injective.
 By Theorem~4, in this case the $B$\+module $H^*(D,C)$ is Koszul
if and only if the $\gr_ND$\+comodule $\gr_NC$ is Koszul.
 The last assertion now follows from Theorem~3, Corollary~5\,(a),
and the version of Corollary~3\,(a) for morphisms of graded
coalgebras. \qed

\begin{thm}{Corollary 9}
 Let $\psi\:G'\rarrow G''$ be a homomorphism of pro-$l$-groups
and $G$ be its kernel.
 Assume that the cohomology algebras $A=H^*(G'',\.\Z/l)$ and
$B=H^*(G',\.\Z/l)$ are Koszul.
 In this case, if the $A$\+module $B$ is Koszul, then
the homomorphism~$\psi$ is injective; if $H_{i,j}(A,B)=0$ 
for all $i-j\ne 0$,~$1$, then the group~$G$ is a free
pro-$l$-group.%
\end{thm}

\pr{Proof}:
 For any morphism of pro-finite groups~$\psi\:G'\rarrow G''$
with the kernel~$G$ the natural morphism of the group coalgebras
$\Z/l(G'')\rarrow \Z/l(G')$ defines the structures of left and
right comodule over the coalgebra $\Z/l(G')$ on the vector
space $\Z/l(G'')$.
 Just as in the proof of Corollary~6, one can see that
the coalgebra $\Z/l(G'')$ is a cofree comodule over $\Z/l(G'/G)$,
the comodule $\Z/l(G'/G)$ over the coalgebra $\Z/l(G')$ is induced
from the trivial $\Z/l(G)$\+comodule, and it follows that
$H^*(\Z/l(G'),\.\Z/l(G''))\simeq
H^*(\Z/l(G'),\.\Z/l(G'/G))\oa_{\Z/l(G'/G)}\Z/l(G'')\simeq
H^*(G,\.\Z/l)\oa_{\Z/l(G'/G)}\Z/l(G'')$.
 So one has $H^i(G,\.\Z/l)=0$ if and only if
$H^i(\Z/l(G'),\.\Z/l(G''))=0$.
 It remains to apply Theorem~8; one does not have to
distinguish $H^{i,j}(A,B)$ from $H^{i,j}(A^\opp,B^\opp)$
because the pro-finite group cohomology algebras are
(skew)\+commutative.  \qed

\begin{thm}{Theorem 9}
 Let $G'$ be a pro-$l$-group, $G\sub G'$ a normal subgroup,
and $G''=G'/G$ the quotient group.
 Denote by $A=H^*(G'',\.\Z/l)$ and $B=H^*(G',\.\Z/l)$ the cohomology
algebras of the quotient group~$G''$ and the group~$G'$.
 Suppose that
  \begin{itemize}
    \item[(1)] the natural homomorphism $f\:A\rarrow B$ is an isomorphism
          in degree~1 and an epimorphism in degree~2;
    \item[(2)] the kernel~$J$ of the morphism~$f$ is isomorphic to
          the shift $K(2)$ of a graded $A$\+module~$K$ for which
          the morphism $r_{A,K}\:\q_AK\rarrow K$ is an isomorphism
          in degree~1 and a monomorphism in degree~2;
    \item[(3)] the algebra~$A$ is Koszul and the $A$\+module $\q_AK$
          is Koszul.
  \end{itemize}
  Then the subgroup~$G$ is a free pro-$l$-group.
\end{thm}

\pr{Proof}:
 Consider the Serre--Hochschild spectral sequence
 $$
   E_2^{p,q}=H^p(G'',H^q(G,\.\Z/l))\implies H^{p+q}(G',\.\Z/l),
   \qquad d_r^{p,q}\:E_r^{p,q}\rarrow E_r^{p+r,q-r+1}.
 $$
 The operator~$d_2\:E_2^{p-2,\,1}\rarrow E_2^{p,0}$ together with
the projection of $E_2^{p,0}$ into $E_\infty^{p,0}$ define a $3$-term
complex of graded modules over $E_2^{*,0}=H^*(G'',\.\Z/l)$
  $$
   H^*(G'',H^1(G,\.\Z/l))\.(2)\lrarrow H^*(G'',\.\Z/l)
   \lrarrow H^*(G',\.\Z/l).
  $$
 Let us denote by $P$ the $G''$\+module $H^1(G,\.\Z/l)$.
 Taking the kernel of the second arrow of this complex, we get
a morphism of graded $A$\+modules $H^*(G'',P)\rarrow K$.

 Since the morphism $H^*(G'',\.\Z/l)\rarrow H^*(G',\.\Z/l)$ is
surjective in the degrees~$1$ and~$2$, one has
$E_3^{0,1}=E_3^{1,1}=0$.
 It follows easily that the map $H^0(G'',P)\rarrow K_0$
is an isomorphism and the map $H^1(G'',P)\rarrow K_1$
is injective.
 Therefore, the quadratic parts of the $A$\+modules $H^*(G'',P)$
and~$K$ coincide.
 Considering the composition
 $$
  \q_AK\,\simeq\,\q_A H^*(G'',P)\lrarrow H^*(G'',P)\lrarrow K
 $$
and using the condition~(2), we deduce that the map
$H^1(G'',P)\rarrow K_1$ is an isomorphism and the morphism
$r_{A,\,H^*(G'',P)}\:\q_A H^*(G'',P)\rarrow H^*(G'',P)$
is an isomorphism in degree~1 and a monomorphism in degree~2.

 According to Theorem~4, it follows from the latter that
the $A$\+module $H^*(G'',P)$ is quadratic, since the algebra~$A$
and the $A$\+module $\q_A H^*(G'',P)\simeq\q_AK$ are Koszul.
 The above composition now shows that the map $H^2(G'',P)\rarrow K_2$
is injective.
 Returning back to the spectral sequence, we conclude that
the morphism $d_2^{0,2}\:E_2^{0,2}\rarrow E_2^{2,1}$ is zero
because the map $H^2(G'',P)\rarrow K_2$ is injective and that
the morphism $d_3^{0,2}\:E_3^{0,2}\rarrow E_3^{3,0}$ is zero
because the map $H^1(G'',P)\rarrow K_1$ is surjective.

 Because the morphism $H^2(G'',\.\Z/l)\rarrow H^2(G',\.\Z/l)$
is surjective, we have $E_4^{0,2}=0$ and it now follows
that $E_2^{0,2}=0$.
 This means that $H^0(G'',H^2(G,\.\Z/l))=0$.
 Since $G''$ is a pro-$l$-group, we have $H^2(G,\.\Z/l)=0$.
 Finally, since the group~$G$ is a pro-$l$-group, it is
a free pro-$l$-group according to a result from~\cite{Ser}.
\qed

\begin{thm}{Lemma 3}
 Let $F$ be a field containing a primitive $l$\+root of unity if $l$
is odd and containing a square root of~$-1$ if\/ $l=2$.
 Then the cohomology algebra of the Galois group $\Gal(\Fradl/F)$ with
constant coefficients $\Z/l$ is naturally isomorphic
to the exterior algebra $\Lambda_l^*(F)=\bwe_{\Z/l}^*(F^*/F^{*l})$.
\end{thm}

\pr{Proof}:
 By Kummer's theory, it is clear that the first cohomology group of
$\Gal(\Fradl/F)$ is isomorphic to $F^*/F^{*l}$.
 The square root of~$-1$ condition guarantees that the cohomology
algebra is skew-commutative.
 It suffices to show that the cohomology algebra is the exterior 
algebra over the first cohomology group.

 Introduce the notation for fields $K=\Fradl$ and $E=\Frootone$.
 By Kummer's theory, the Galois group $\Gal(K/E)$ is naturally
a subgroup of the group of homomorphisms from~$F^*$ to the projective
limit of the groops of $l^N$\+roots of unity, $\varprojlim\mu_{l^N}$.
 The group $\Hom(F^*,\,\varprojlim\mu_{l^N})$ is an abelian pro-$l$-group
without torsion, hence so is the group $\Gal(K/E)$.
 It follows that $\Gal(K/E)$ is a free abelian pro-$l$-group and its
cohomology ring is an exterior algebra.

 If the field $F$ contains all the $l^N$\+roots of unity, then $E=F$ and
we are done.
 Otherwise, the group $\Gal(E/F)$ is isomorphic to $\Z_l$ and it acts
trivially in the cohomology algebra of $\Gal(K/E)$.
 The Serre--Hochschild spectral sequence degenerates into a series of
exact sequences of the form
$$
 0\rightarrow H^{n-1}(\Gal(K/E),\.\Z/l) \rarrow H^n(\Gal(K/F),\.\Z/l) 
 \rarrow H^n(\Gal(K/E),\.\Z/l)\rightarrow0.
$$
It follows easily that $H^*(\Gal(K/F),\.\Z/l)$ is an exterior algebra.
\qed

\smallskip

\pr{Proof of Theorem 2}:
 It follows from the condition~(1) of Theorem~2 that the
quadratic part $\q H^*(G_F,\.\Z/l)$ of the Galois cohomology
algebra is isomorphic to the Milnor algebra $\KM_*(F)\ot\Z/l$.
 By Theorem~6, it follows from the condition~(2) that
the latter algebra is Koszul.
 Therefore, Corollary~2 yields the isomorphisms
  $$
   H^*(G_F^{(l)},\.\Z/l)\simeq \q H^*(G_F,\.\Z/l)
   \simeq \KM_*(F)\ot\Z/l.
  $$

 Now it remains to apply Theorem~9 or Corollary~9
to the pro-$l$-group $G'=G_F^{(l)}$ with the quotient group
$G''=\Gal(K/F)$ and the corresponding subgroup
$G=\Gal(F^{(l)}/K)$ (where $F^{(l)}$ is the maximal
Galois pro-$l$-extension of~$F$).
It is obvious that $G$ is isomorphic to the maximal
quotient pro-$l$-group of the absolute Galois group
of~$K$ that we are interested in.
\qed

\appendix

\bigskip
\section*{Appendix. Supporting Evidence for Conjecture 1}
\medskip

 In this section we present a counterexample showing that the most
straightforward generalization of Bogomolov's conjecture to
arbitrary fields fails and prove several results providing evidence
in support of our Conjecture~1.

\begin{thm}{Counterexample}
 There exist fields~$F$ such that the commutator subgroup of the
maximal quotient pro-$l$-group of the absolute Galois group
of~$F$ is not a free pro-$l$-group.
 Moreover, there are fields such that the commutator subgroup
of their Sylow pro-$l$-subgroup of the absolute Galois group
is not free.
\end{thm}

\pr{Proof}:
Let $F$ be a field of characteristic not equal to~$l$ which
contains a primitive root of unity of degree~$l$ but no primitive
roots of unity of degree~$l^N$ for some $N>1$ 
(assumed to be the minimal such integer).
 Suppose that the maximal quotient pro-$l$-group of~$G_F$
is not abelian.
 It is very easy to find a number field, $l$-adic field, or
power series field satisfying the listed properties.

 Consider the field $F((z))$ of formal Laurent power series in
the variable~$z$ with coefficients in~$F$.
 Then the maximal quotient pro-$l$-group of the absolute Galois
group of the field of power series $G_{F((z))}^{(l)}$ is isomorphic
to the semidirect product of the group $G_F^{(l)}$ with the group
of $l$-adic integers $\Z_l$, where $G_F^{(l)}$ acts via the cyclotomic
character.
 The commutator subgroup of this semidirect product group is
the semidirect product of the commutator subgroup of $G_F^{(l)}$
with $l^N\Z_l$.
 This is clearly not a free pro-$l$-group.
 The Sylow pro-$l$-subgroup of the absolute Galois group of $F((z))$
is isomorphic to the semidirect product of the Sylow subgroup of $G_F$
with $\Z_l$.
 The commutator subgroup of this group is described as above;
it is not a free pro-$l$-group either.

 In particular, for any field $F$ containing a primitive $l$\+root
of unity but no primitive $l^N$\+roots of unity for some $N>1$,
the field of doubly iterated formal Laurent power series $F((w))((z))$
provides the desired counterexample.
\qed

\begin{thm}{Proposition A.1}
 Conjecture~1 holds for number fields and their functional analogues,
i.~e., algebraic extensions of\/~$\Q$ and\/ $\F_q(x)$.
 Moreover, in these cases it suffices to add roots of unity only, i.~e.,
take $K=\Frootone$.
\end{thm}

\pr{Proof}:
 The assertion is (essentially) a particular case of the result
of~\cite[II.3.3]{Ser}.
 It suffices to prove that $H^2(G_M,\.\Z/l)=0$ for any field $M$ between
$K$ and $\oF$.
 For any field $L$ containing a primitive $l$-root of unity, the group
$H^2(G_L,\.\Z/l)$ is isomorphic to the subgroup of the Brauer group
$\Br(L)$ annihilated by~$l$.
 Obviously, the Brauer group functor preserves direct limits of fields.
 So it suffices to show that for any field $E$ finite over~$\Q$ or
$\F_q(x)$ and any element $\alpha\in\Br(E)$ there exists a field $L$
finite over~$E$ and contained in the composite $EK$ such that the image
of~$\alpha$ in $\Br(L)$ vanishes.

 It remains to use the description of Brauer groups given by class field
theory.
 The Brauer group of a global field is a subgroup of the direct sum of
the Brauer groups of local completions.
 For any class $\alpha\in\Br(E)$ there is a finite number of completions
$E_v$ of the field~$E$ where the local class $\alpha_v$ is nonzero.
 Adjoining $l^N$\+roots of unity with large enough~$N$ to a nonarchimedean
complete field $E_v$, one can obtain extensions of $E_v$ of degrees
divisible by arbitrary high powers of~$l$ (since $E_v$ is finite over
a field $\Q_p$ or $\F_q[[z]]$, which has this property).
 And by local class field theory an $l$-torsion element of $\Br(E_v)$ dies
in any finite extension of~$E_v$ of degree divisible by~$l$.
 In the archimedean case, the local Brauer group is zero (for $E_v=\C$)
or $\Z/2$ (for $E_v=\R$), and it suffices to adjoin the square root
of~$-1$ in order to kill it.
\qed

\begin{thm}{Proposition A.2}
 Let $F$ be a Henselian discrete valuation field and $f$~be its residue
field.
 Then the field $F$ satisfies Conjecture~1 whenever $f$ does
(for a given~$l$).
 In particular, if\/ $\chr f=l$, then $F$ satisfies Conjecture~1.
\end{thm}

\pr{Proof}:
 As in the proof of Proposition~A.1, we will argue in terms of the Brauer
groups.
 The Brauer groups of Henselian discrete valuation fields with perfect
residue fields were computed in \cite[chapitre XII]{Ser2}
and~\cite[\S2]{Pop}.

 Consider the two cases $\chr f\ne l$ and $\chr f=l$ separately.
 It is easy to verify following the arguments of \cite{Ser2}
and~\cite{Pop} that for $\chr f\ne l$ there is always a natural exact
sequence
$$
 0\lrarrow \Br(f)_{(l)} \lrarrow \Br(F)_{(l)} \lrarrow 
 \Hom(G_f,\.v(F)\otimes\Q/\Z)_{(l)} \lrarrow 0,
$$
where $v(F)$ denotes the group of values of the discrete valuation,
$G_f$ is the absolute Galois group of~$f$, and the subindex~$l$ in
parentheses signifies the $l$-primary components of the corresponding
abelian groups.
 Let $E$ be a finite extension of the field $F$ and $e$ be its
residue field.
 We would like to prove that any element of $\Br(E)_{(l)}$ can be killed
by a radical extension of~$F$.
 Let $\alpha\in\Br(E)_{(l)}$; for compactness reasons, the corresponding
element of $\Hom(G_e,\.v(E)\otimes\Q/\Z)$ is a homomorphism with finite
image.
 We kill this class by adjoining a root of an appropriate degree from
a uniformizing element of the valuation~$v$ in~$F$.
 This extends our field $E$ to a larger field $D$, still finite over $F$;
let $d$ denote its residue field.
 Applying the above exact sequence to the valued field $D$, we get
an element $\alpha'\in\Br(d)_{(l)}$.
 By our assumption, this element can be killed by adjoining $l$-power
roots of elements of~$f$.
 To make the original element $\alpha\in\Br(E)_{(l)}$ vanish, it now
suffices to adjoin the corresponding roots of any preimages in~$F$ of
those residue elements, together with the root of the uniformizing
element in $F$ that we adjoined already.

 In the case $\chr f = l$ the main step is to reduce the problem to
the situation of a perfect residue field.
 Let us repeat transfinitely the following process of constructing purely
inseparable extensions~$e$ of the field~$f$.
 Choose any element of~$f$ which is not an $l$\+power in~$f$ and adjoin
to $f$ all roots of this element of degrees~$l^N$; this is the first step.
 On each of the subsequent steps we choose an element of~$f$ which is
still not an $l$-power in~$e$ and adjoin all its $l^N$\+degree roots to~$e$.
 When the transfinite induction terminates, we have a purely inseparable
extension~$d$ of the field~$f$ which is easily found to be its perfect
closure.
 Indeed, any $x\in f$ is an $l$\+power in~$d$ by construction;
let us assume that it is an $l^{n-1}$\+power and show that it is also
an $l^n$\+power.
 Let $x=y^{l^{n-1}}$, \ $y\in d$.
 The element~$y$ can be expressed in terms of elements of~$f$ and roots
of elements of~$f$ that we adjoined in the process.
 But both kinds of elements are $l$\+powers in~$d$ by the construction.

 Now we repeat this transfinite process on the level of valued fields.
 At each step we choose a compatible family of roots of all degrees $l^N$,
one root of each degree, taken from an arbitrary element of $F$ 
lifting $x\in f$.  
 Notice that we avoid any extension of the group of values of
the valuation; indeed, at each step the total degree of the extention
coincides with the degree of extension of the residue fields.
 So the valuation remains discrete.
 The transfinite process results in the construction of a Henselian
discrete valuation field~$D$ with a perfect residue field.
 According to \cite{Ser2} and~\cite{Pop}, we have again an exact sequence
$$
 0\lrarrow \Br(d) \lrarrow \Br(D) \lrarrow 
 \Hom(G_d,\.v(D)\otimes\Q/\Z) \lrarrow 0.
$$
 In addition, the $l$\+primary component of the Brauer group of a perfect
field of characteristic~$l$ vanishes.
 By construction, we have $F\sub D\sub\Fradl$.

 Finally, let $M$ be a finite extension of the field $D$ and $m$ be its
residue field.
 It remains to annihilate the class from
$\Hom(G_m,\.v(M)\otimes\Q/\Z)_{(l)}$.
 As above, we do this by adjoining a root of high enough $l$\+power
degree from a uniformizing element of the valuation~$v$ of the field~$F$.
\qed

\begin{rem}{Remark 1:}
 Here is another partial result in support of Conjecture~1, provided by
a version of an argument which I learned from F.~Bogomolov sometime
in 1995\+-96.
 Let $F$ be a field of characteristic not equal to~$l$; set
$K=\Fradl$.
 The conjecture can be rephrased by saying that for any finite extension
$E$ of the field $F$ the image of the group $\KM_2(E)/l$ in
$\KM_2(KE)/l$ vanishes.
 Let us show that this is true for any finite Galois extension $E/F$
with the Galois group isomorphic to a direct sum of copies of $\Z/2$,
or a direct sum of copies of $\Z/3$, or, more generally, for any
abelian field extension with the Galois group annihilated by a number~$q$
such that the field $F$ has no nontrivial finite extensions of
degrees less or equal to~$q/2$.
 In particular, the assertion applies to any extension $E/F$ with
the Galois group $\Gal(E/F)\simeq \Z/l\op\dsb\op\Z/l$ of a field $F$
having no finite extensions of the degrees prime to~$l$.

 The map $\KM_2(KE)/l\rarrow\KM_2(KE)\ot_\Z\Q_l/\Z_l$ being injective,
it suffices to prove that the image of the group $\KM_2(E)\ot_\Z\Q_l/\Z_l$
in $\KM_2(KE)\ot_\Z\Q_l/\Z_l$ vanishes.
 Now $\KM_2(E)\ot_\Z\Q$ is a representation of the finite group $\Gal(E/F)$
over a field of characteristic zero; so it decomposes into a direct sum
of irreducibles.
 Moreover, any element of $\KM_2(E)\ot_\Z\Q$ invariant under the action of
a subgroup of $\Gal(E/F)$ comes from an element of $\KM_2(L)\ot_\Z\Q$
for the corresponding subfield $F\sub L\sub E$ (which can be obtained
using the transfer operation).
 If the Galois group $\Gal(E/F)$ is abelian, its irreducible
representations are actually representations of its cyclic quotient
groups, so their elements come from cyclic extensions $L/F$ of
degrees dividing~$q$.
 
 By the Bass--Tate lemma, for any extension $L/F$ of degree~$\le q$
of a field~$F$ having no nontrivial finite extensions of
degrees~$\le q/2$, the multiplication map $\KM_1(L)\ot_\Z\KM_1(F)
\rarrow\KM_2(L)$ is surjective~\cite[section~4]{Pos}.
 Hence the map $\KM_1(L)\ot_\Z\KM_1(F)\ot_\Z\Q_l/\Z_l\rarrow
\KM_2(L)\ot_\Z\Q_l/\Z_l$ is surjective, too.
 Given that the map $\KM_1(F)\ot_\Z\Q_l/\Z_l\rarrow
\KM_1(K)\ot_\Z\Q_l/\Z_l$ is zero by definition, it follows that the map
$\KM_2(L)\ot_\Z\Q_l/\Z_l\rarrow\KM_2(KL)\ot_\Z\Q_l/\Z_l$ induced by
the embedding of fields $L\rarrow KL$ also vanishes.

 Summing up the two previous paragraphs, since the group
$\KM_2(E)\ot_\Z\Q$ is generated by the images of the groups
$\KM_2(L)\ot_\Z\Q$ for intermediate fields $F\sub L\sub E$
of degree~$\le q$ over~$F$ and the map $\KM_2(E)\ot_\Z\Q\rarrow
\KM_2(E)\ot_\Z\Q_l/\Z_l$ is obviously surjective, we can conclude that
the map $\KM_2(E)\ot_\Z\Q_l/\Z_l\rarrow\KM_2(KE)\ot_\Z\Q_l/\Z_l$ is zero.
 The same argument also applies, e.~g., to nonabelian extensions $E/F$
of any field $F$ with the Galois group $\Gal(E/F)$ isomorphic to
the symmetric group~$\mathbb S_3$, as its only irreducible representation
of dimension~$>1$ over~$\Q$ is generated by its three 1\+dimensional
subspaces of invariants of the transpositions in~$\mathbb S_3$, which
correspond to intermediate subfields $L$ of degree~$3$ over~$F$.
\end{rem}

\begin{rem}{Remark 2:}
 The following assertion describing the absolute Galois groups of fields
as being ``close to free profinite groups'' in a sense quite different
from that of Conjectures~1\+-2 is not difficult to prove.
 Let $F$ be a finitely generated field and $G=\Gal(\oF/F)$ be its
absolute Galois group.
 Then the group $G$ admits a finite decreasing filtration $G=G^{-1}
\supset G^0\supset G^1\supset\dsb\supset G^N\supset G^{N+1}=\{1\}$ by
closed subgroups normal in $G$ such that the quotient group $G/G^0$
is isomorphic to either $\{1\}$ or $\Z/2$, while all the subsequent
quotient groups $G^n/G^{n+1}$, \ $n\ge0$ are closed subgroups of free
pro-finite groups.
 Moreover, the group $G^0/G^1$ can be only nontrivial when the field
$F$ has characteristic~$0$ and does not contain the square root of~$-1$.
 Replacing finite decreasing filtrations of profinite groups with
transfinitely descreasing ones, one extends the result to arbitrary
(infinitely generated) fields~$F$.

 It suffices to consider the case of a purely transcendental field
$F=\F_p(x_1,\dotsc,x_N)$ or $\Q(x_2,\dotsc,x_N)$.
 In the finite characteristic case, let us filter the field $F$ by its
subfields $F_n=\F_p(x_1,\dotsc,x_{n-1})$ and the field $\oF$ by its
subfields $L_n=\oF_nF$.
 The fields $F=L_0\sub L_1\sub L_2\sub\dsb\sub L_N\sub L_{N+1}=\oF$
are intermediate fields between $F$ and $\oF$, normal over~$F$.
 Set $G_n=\Gal(\oF/L_n)\sub G=\Gal(\oF/F)$ for $n=0$,~\dots, $N+1$.
 The Galois group $G^0/G^1=\Gal(L_1/F)\simeq\Gal(\overline\F_p/\F_p)=
\widehat\Z$ is a free pro-finite group with one generator,
and the Galois groups $G^n/G^{n+1}=\Gal(L_{n+1}/L_n)\simeq
\Gal\big(\overline{\oF_n(x_n)}/\oF_n(x_n)\big)$ for $n\ge1$
are known to be free pro-finite groups~\cite[Corollary~11.34]{Vol}.
 
 In the zero characteristic case, let us first filter the field
$\overline\Q$ by its subfields $\Q=\overline\Q_{-1}\sub\Q_0\sub\Q_1
\sub\Q_2=\overline\Q$.
 Set $\Q_0=\Q[\sqrt{-1}]$. 
 To define the field $\Q_1$, consider the maximal abelian extension
$\Q^\ab=\Q[\sqrt[\infty]{1}]$.
 The Galois group $\Gal(\Q^\ab/\Q)$ is naturally isomorphic to
the product over prime numbers $\widehat\Z^*=\prod_l\Z_l^*$ of
the multiplicative groups of $l$\+adic integers.
 Its maximal torsion subgroup is $\widehat\Z^*_{\mathrm{tors}}\simeq
\Z/2\times\prod_l(\Z/l)^*$.
 Consider the intermediate field $\Q_1'$ between $\Q$ and $\Q^\ab$
corresponding to the closed subgroup $\widehat\Z^*_{\mathrm{tors}}
\subset\widehat\Z^*$, and set $\Q_1=\Q_0\Q_1'$.
 Then one has $\Gal(\Q_0/\Q)\simeq\Z/2$ and $\Gal(\Q_1/\Q_0)\simeq
\widehat\Z$.
 The pro-finite group $\Gal(\overline\Q/\Q_1)$, having cohomological
dimension~$1$ by a slight strengthening of the result of
Proposition~A.1 above (provable by the same method), is
a closed subgroup of a free pro-finite group~\cite[section~22.4]{FJ}.

 Finally, we filter the field $F$ by its subfields $F_n=
\Q(x_2,\dotsc,x_{n-1})$, \ $n=2$,~\dots, $N+1$, and the field $\oF$ by
its subfields $L_0=\Q_0F$, \ $L_1=\Q_1F$, and $L_n=\oF_nF$ for $n\ge2$.
 The fields $F=L_{-1}\sub L_0\sub L_1\sub\dsb\sub L_N\sub L_{N+1}=\oF$
are intermediate fields between $F$ and $\oF$, normal over~$F$.
 Set $G_n=\Gal(\oF/L_n)\sub G=\Gal(\oF/F)$ for $n=-1$,~\dots, $N+1$.
 Then the group $G^{-1}/G^0=\Gal(L_0/F)\simeq\Gal(\Q_0/\Q)$ is isomorphic
to $\Z/2$, the group $G^0/G^1=\Gal(L_1/L_0)\simeq\Gal(\Q_1/\Q_0)$ is
a free pro-finite group with one generator, the group
$G^1/G^2=\Gal(L_2/L_1)$ is isomorphic to $\Gal(\overline\Q/\Q_1)$,
and the groups $G^n/G^{n+1}=\Gal(L_{n+1}/L_n)\simeq
\Gal\big(\overline{\oF_n(x_n)}/\oF_n(x_n)\big)$ for $n\ge2$
are free pro-finite groups.
\end{rem}

\end{document}